\documentclass{amsart}

\usepackage[T1]{fontenc}
\usepackage{graphicx}
\usepackage{hyperref}

\title{Data, geometry and homology} 

\author{Jens Agerberg}
\address{KTH, Department of Mathematics, S 10044 Stockholm, Sweden}
\email{jensag@kth.se}

\author{Wojciech Chach\'olski}
\address{KTH, Department of Mathematics, S 10044 Stockholm, Sweden}
\email{wojtek@kth.se}

\author{Ryan Ramanujam}
\address{KI, Department of Clinical Neuroscience, S 17176 Stockholm, Sweden}
\email{ryan.ramanujam@ki.se}

\begin{document}

\maketitle

\begin{abstract}
Homology-based invariants can be used to characterize the geometry of datasets and thereby 
 gain some understanding of 
the  processes generating  those datasets. In this work we investigate how the geometry of a dataset changes when it is subsampled in various ways. In our framework the dataset serves as a reference object; we then consider different points in the ambient space and endow them with a geometry defined in relation to the reference object, for instance by subsampling the dataset proportionally to the distance between its elements and the point under consideration. We illustrate how this process can be used to extract rich geometrical information, allowing for example to classify points coming from different data distributions.
\end{abstract}

\section{Introduction}
\label{sec:Introduction}
Consider a point in the plane. In itself a point does  not have an interesting geometry, however in relation to  other objects (called  {\em reference} objects) it has   rich  geometrical aspects
such as being on  the left or right side of an oriented line, or being inside or outside a circle.  These aspects can of course become quite complicated and difficult to interpret when the reference object is complex.
The aim of this note is to present a strategy of how  such relative geometrical aspects of points can be effectively
encoded by  homology-based invariants and illustrate how these invariants  can    be used to distinguish different types of points. For example, consider the very classical problem of classifying handwritten digits in MNIST data~\cite{MNIST}. To decide if a  handwritten digit represents $1$ or $7$,
we might look at the geometrical aspects of the point representing the digit relative to, for instance,
the reference object  formed by some representatives of other handwritten digits. Standard methods provide a very high classification accuracy for MNIST data. Our purpose is to illustrate  cases where homology-based invariants describing  some of these  relative geometrical aspects of points representing handwritten digits also contain a large amount of discriminative information. 
For both reading and writing purposes, humans have learned which  variations of written digits can be recognized. Thus the key information  needed to  identify a handwritten digit is encoded in the spaces  of points representing each digit.   It is therefore not  surprising that    geometrical aspects relative to these spaces   provide accurate classification, which is well illustrated by for example high precision of k-nearest neighbors (KNN)  method. However, it is surprising  that {\em homological} invariants extracted from  relative geometrical information also contain a large amount of discriminatory information. 

What do we mean by geometrical information?
In this article  it is information extracted from any space associated with the considered problem or object, where  by a space we mean a simplicial complex or a family of them. 
For example let $\mathcal{R}$ be a finite subset of $\mathbb{R}^n$.  By restricting a metric from $\mathbb{R}^n$ to $\mathcal{R}$  we can form  its Vietoris-Rips complexes~\cite{MR1368659}. 
This is just one instance of a possible spacial representation of metric properties of  $\mathcal{R}$. 
Here is another possibility which we explore in this article. Choose a natural number $s$ (called {\em sample size}) and consider the collection  of the Vietoris-Rips complexes of some  $s$-element subsets of $\mathcal{R}$, chosen according to some specified rule, for instance selecting 
all of them, or only those that are at certain distance to a given point. 
We find that this type of spacial representation of metric properties of $\mathcal{R}$ is  informative.
The way we use this representation is as follows. 
First, persistent homology is extracted  for each of the subsets (there is an extensive literature regarding persistent homology see for example~\cite{ MR2405684, MR2358377,  MR2777589} and their bibliographies). 
If $s$ is relatively small, the computational cost of this step is reasonable and can be done in parallel for each subset and for a rather high homological degree. The next step is to average these outcomes over some choice of the $s$-element subsets.  This step requires transforming persistent homologies into objects whose averages (expected values) can be calculated. For that purpose we utilise  stable ranks~\cite{MR4057607, OliverWojtek, MR3735858} that represent persistence modules by non-increasing piecewise  constant functions (see Appendix~\ref{sec:adfgdfagsdfg}).  Stable ranks also enable us to  use the associated  kernel~\cite{10.3389/fams.2021.668046} and for example SVMs for classification purposes.
We should mention that there are other ways of representing persistence modules by objects suitable for statistical analysis, for example persistence landscapes~\cite{MR3317230} or persistence images~\cite{MR3625712}. See also~\cite{chazal:hal-01073073} where subsampling techniques in a similar context have been explored.

 The aim of this article is   to present a pipeline for assigning  a stable rank (non-increasing 
  piecewise  constant function) to a subset $\mathcal{R}$ in $\mathbb{R}^n$.  The space of parameters for our pipeline naturally  splits into two types  called {\em global} and {\em relative}. Global stable ranks encode some geometrical aspects of  $\mathcal{R}$. Relative stable ranks encode some geometrical aspects of points in  $\mathbb{R}^n$ relative to $\mathcal{R}$. 
   In Section~\ref{sec:asgwdfhfg} we give simple examples in which it is possible to geometrically interpret the outcomes of the relative pipeline (lying on specified side of a hyperplane or inside  a circle).
   Although in general geometrical interpretation of
   both relative and global  stable ranks  may be too complex, these invariants can be used for distinguishing purposes.  For example we show that
   the training set in MNIST representing  $7$ is geometrically different from the test set in MNIST representing  $7$, as their invariants based on the  homology in degree 1 are quite different. Thus, the geometry of the testing dataset for $7$ is not entirely representative of the geometry of its training dataset.
   In Section~\ref{sec:global} we also  illustrate variability among  global stable ranks of MNIST training datasets    across different digits  in homological degrees $0$, $1$, and $2$. We explore this variability in Section~\ref{sec:asgwdfhfg}
   for classification purposes. For example if as a reference object we choose the union of the training MNIST datasets for digits $1$ and $7$, then the test digits labeled by $1$ and $7$ can be quite accurately classified using
   their relative stable ranks. We believe this is a consequence of the fact that global geometrical aspects of the  training MNIST datasets for $1$ and $7$, as measured by their global stable ranks, are quite different.  By repeating analogous experiments for other pair of digits, when global stable ranks of the digits are similar, it becomes harder to distinguish the test digits by using
   analogous  parameters and reference objects.

Extracting stable ranks is a simplifying procedure whereby a large amount of information is  discarded. The challenge is to be able to steer some of the choices of the parameters which control stable ranks   in such a way that  some of the aspects relevant to the problem at hand  are retained. In this paper we showcase that this strategy is viable and indeed stable ranks can contain interesting information particularly  if  an appropriate reference object is chosen.

\section{Pipeline}
\label{sec:pipeline}
The initial input  is  a finite subset 
$\mathcal{R}\subset \mathbb{R}^r$ called
a  {\em reference object}. We are going to explain how to assign to it various non-increasing piecewise constant  functions. These functional representations of $\mathcal{R}$ are then used as inputs for various  analysis pipelines
such as SVMs. 

\subsection*{Step A:  probabilities.}

The objective is to obtain a function $\text{prob}\colon \mathcal{R}\to \mathbb{R}$, called {\em probability}, with the following properties: all its values are non-negative, and their sum  $\Sigma_{x\in\mathcal{R}} \text{prob}(x)$ is either $1$ or $0$. For example we could take the uniform probability which is the constant function
with value  $1/|\mathcal{R}|$. Here is another
construction, divided into two steps A1 and A2, of probability functions that play  important role in this article:

\subsection*{Step A1:  filter function.}

The objective is to obtain a function
$\text{filter}\colon \mathcal{R}\to \mathbb{R}$ called
a {\em filter}. In our particular construction of a filter, the input consists of a point $p$ in $\mathbb{R}^r$ and a vector field on $\mathcal{R}$ 
 represented by a function $\mathcal{V}\colon \mathcal{R}\to \mathbb{R}^r$.  In this article the focus is on   two types of a vector field:
 a constant vector field, and a vector field $\mathcal{V}_c\colon \mathcal{R}\to \mathbb{R}^r$ determined by a point $c$, called {\em center}, in $\mathbb{R}^r$  which assigns to
 $x$ in $\mathcal{R}$ the vector 
 $\mathcal{V}_c(x):=c-x$ from $x$ to $c$.
 For example we could take $c$ to be the point $p$ or the center of mass of $\mathcal{R}$.
 
 In this step the output is the following filter, determined by the point $p$ and the vector field $\mathcal{V}$. For $x$ in $\mathcal{R}$:
\[\text{filter}(x):=\begin{cases}
0 &\text{ if } \mathcal{V}(x)=0\\
|\text{proj}_{\text{span}(\mathcal{V}(x))}(p-x)|&\text{ if }  \mathcal{V}(x)\cdot (p-x )  \geq 0\\
-|\text{proj}_{\text{span}(\mathcal{V}(x))}(p-x)|&\text{ if }  \mathcal{V}(x)\cdot (p-x ) <0
\end{cases}\]
For example, if the consider  vector field is given by $\mathcal{V}_p$, the associated filter function assigns to $x$ in $\mathcal{R}$
the distance between $x$ and $p$.

\subsection*{Step A2: distribution and  probabilities.}\label{sec:pipeline:a2}

In this note a {\em distribution} is 
a function $\mathcal{D}\colon \mathbb{R}\to \mathbb{R}$
whose  values are  non-negative. 

In this step a distribution $\mathcal{D}$ needs to be chosen. It is used to obtain the following  
 probability function  which is the outcome of this step.
Let $S=\Sigma_{y\in \mathcal{R}}\mathcal{D}(\text{filter}(y))$.
\[\text{prob}(x):=
\begin{cases}
0 & \text{ if } S=0\\
\mathcal{D}(\text{filter}(x))/S &
\text{ if }S\not =0
\end{cases}\]
 
\subsection*{Step B: averaged stable ranks.}
 
The objective is to obtain a  non-increasing piecewise constant  function representing the reference object $\mathcal{R}$. 
 
\subsection*{Step B1: sub-sampling.}
 
The  probability function $\text{prob}\colon \mathcal{R}\to \mathbb{R}$, obtained in step A, is  used to 
sample the reference object $\mathcal{R}$.
For this purpose two natural numbers  $s$ and $n$ need to be chosen,   called respectively 
{\em sample size} and {\em number of instances}. Depending on these numbers, the outcome of this step is a set $\mathcal{S}$ described as follows:
\begin{itemize}
    \item If $s>|x\in\mathcal{R}\ |\ \text{prob}(x)>0|$, then the outcome $\mathcal{S}$ is the empty set.
    \item If $s\leq |x\in\mathcal{R}\ |\ \text{prob}(x)>0|$, then the outcome $\mathcal{S}$ is of size $n$ whose elements are subsets of $\mathcal{R}$ of size $s$.
    Each of these subsets is a random choice
    (without replacement) of $s$ elements from $\mathcal{R}$ according to the probabilities specified by the function prob.
\end{itemize}

\subsection*{Step B2: stable ranks.}\label{sec:pipeline:b2}

In this step four choices need to be made.
First, is a natural number $l$ called {\em homological degree}. Second a finite field of coefficient $\mathbb{F}$. Third, a function $f\colon [0,\infty)\to (0,\infty)$ defined on non-negative reals  whose values are strictly positive, called {\em density}. Fourth,   an extended non-negative real number $T$ in $[0,\infty]$, called {\em truncation}.  
These choices are used as follows:
\begin{itemize}
    \item 
Every element $\sigma$ of the outcome $\mathcal{S}$ of  step B1,  which is a subset of the reference object $\mathcal{R}$, is converted  into the following persistence module (the $l$-th homology with coefficients in $\mathbb{F}$ of the corresponding Vietoris-Rips complex, with respect to the Euclidean distance):
\[t\mapsto H_l(\text{VR}_t(\sigma),\mathbb{F})\] 
\item For every $\sigma$ in $\mathcal{S}$, the obtained persistence module is transformed into 
a  non-increasing piecewise constant  function
given by  its {\em stable rank} $\widehat{\text{rank}}(\sigma)$ with respect to the distance type contours $D_f/T$, associated with the density $f$, and  truncated at $T$~\cite[Definitions 5.4 and 5.6]{MR4057607}
(see also Appendix~\ref{sec:adfgdfagsdfg}). 
\item The final outcome of the entire pipeline is the average of all these stable ranks across all $\sigma$ in $\mathcal{S}$:
\[
\widehat{\text{rank}}_{\text{prob},s,n,l,\mathbb{F},f,T}\mathcal{R}:=\left (\Sigma_{\sigma\in \mathcal{S}}\widehat{\text{rank}}(\sigma) \right)/n
\]
\end{itemize}
Throughout the article we use only  $\mathbb{F}_2$ as the field of coefficients since we have found no difference by considering other fields in the presented examples. We also focus on the standard density function given by the constant function $1$.

\section{Global stable ranks}\label{sec:global}

The results of the pipeline described in Section~\ref{sec:pipeline}, when the outcome of step A is given by the uniform probability function, are called {\em global stable ranks} of the reference object. These global stable ranks
encode   aspects of the geometry of the reference object captured by homologies of its $s$-element subspaces. 
In this section we illustrate examples of global stable ranks for the MNIST dataset~\cite{MNIST}.
Recall that  MNIST is a dataset of handwritten digits widely used in machine learning, composed of 60000 training samples and 10000 test samples. The samples are considered as points in $\mathbb{R}^{784}$, since the images have $28 \times 28 = 784$ pixels.
For every $d$ in $\{0,1,\ldots,9\}$, consider two reference objects 
$\text{Test}_d\subset \mathbb{R}^{784}$ and $\text{Train}_d\subset \mathbb{R}^{784}$ formed by these 
handwritten digits in respectively the test  and the training sets of MNIST which are labeled by $d$. 

\begin{figure}[h]
\begin{center}
\includegraphics[height=4cm]{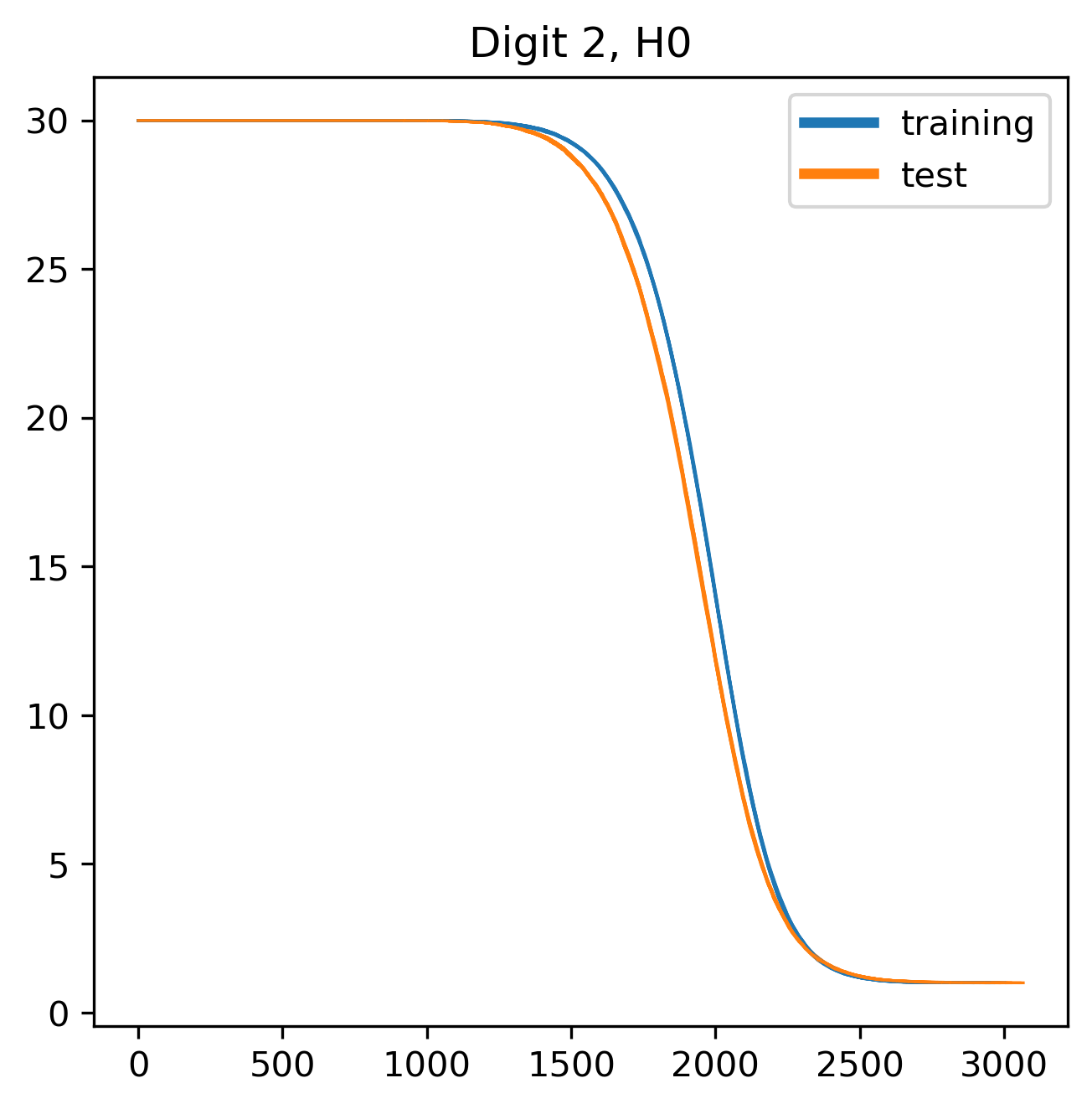}
\includegraphics[height=4cm]{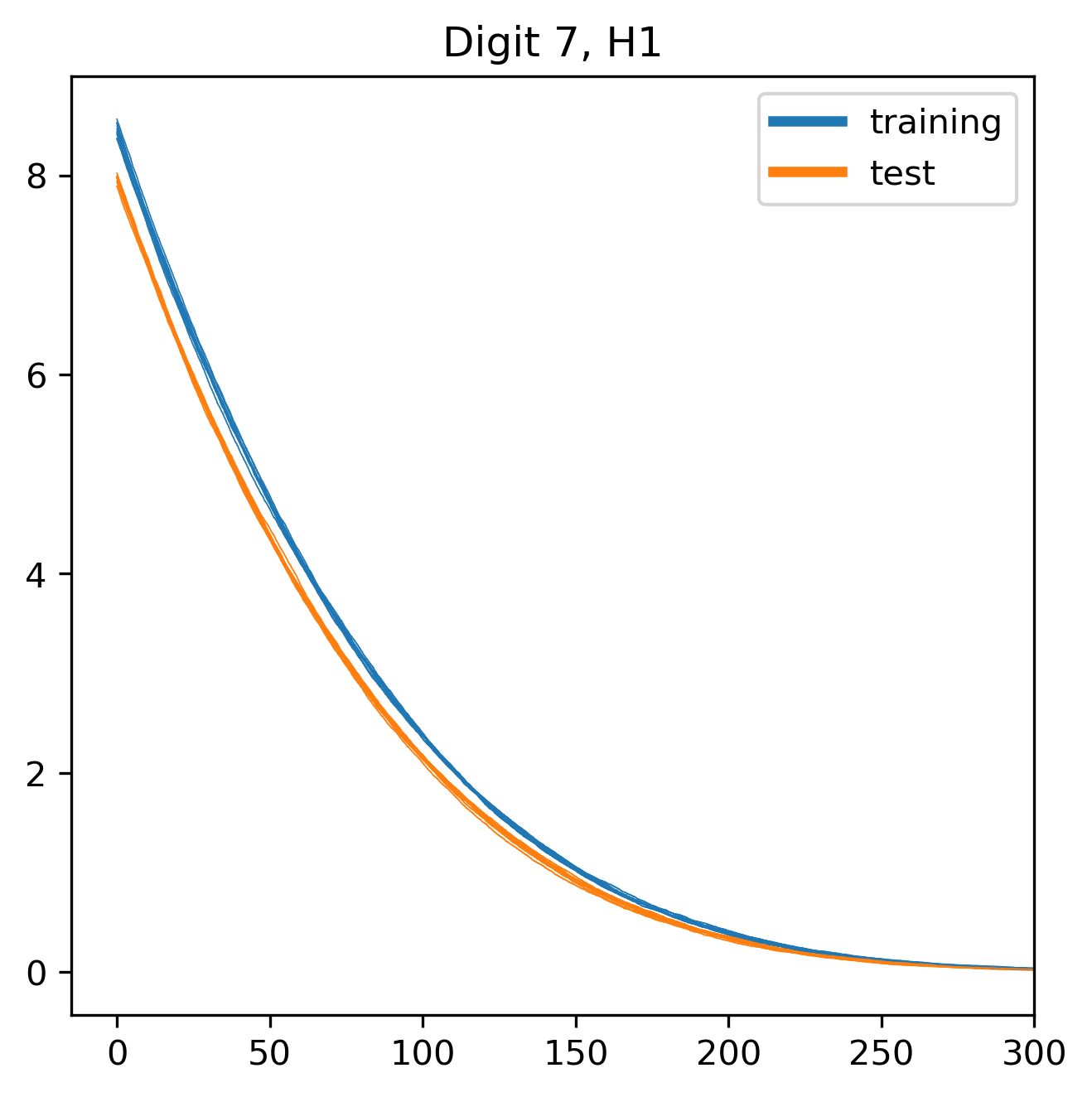}
\includegraphics[height=4cm]{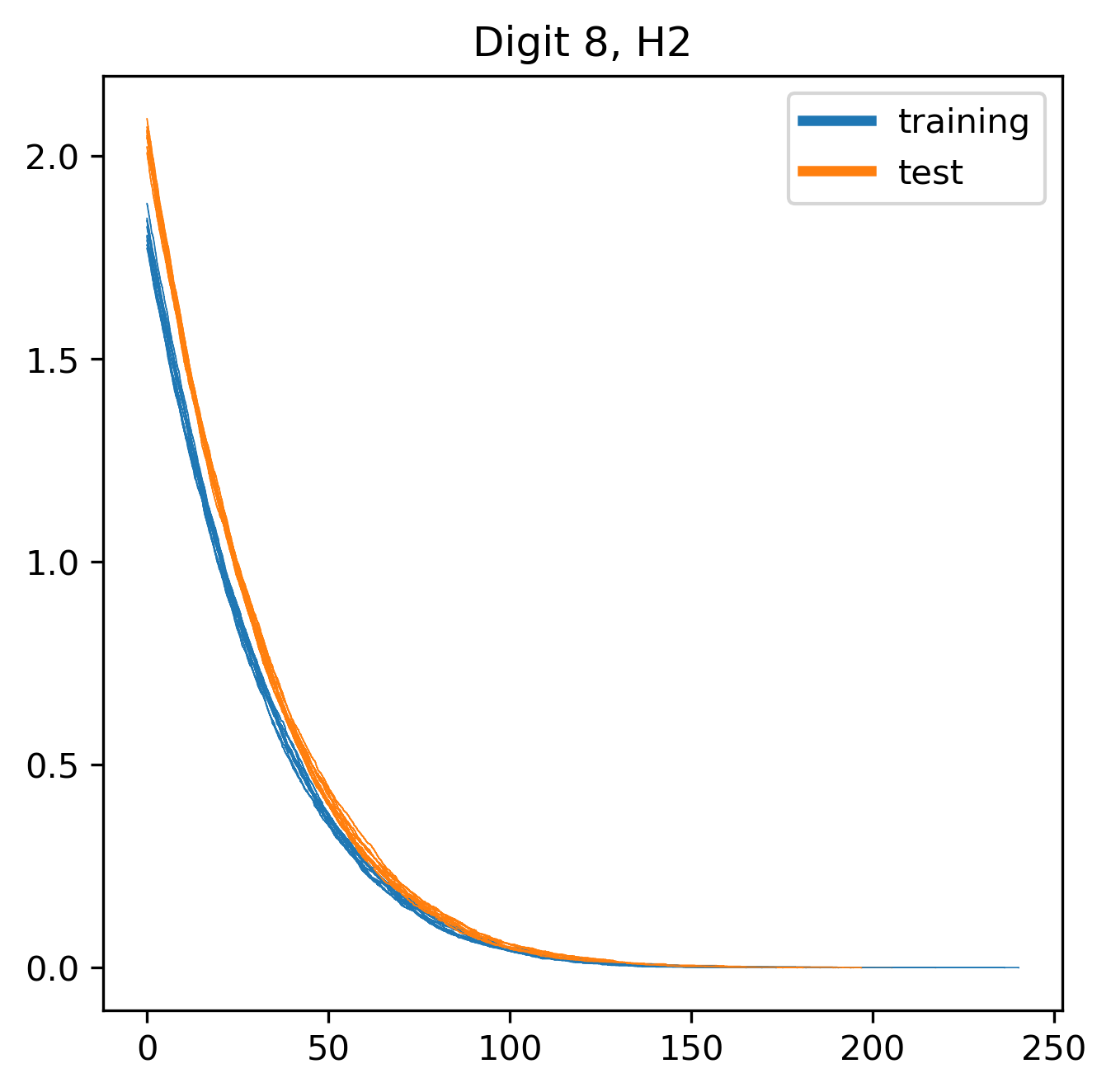}
\end{center}
\caption{The pipeline was repeated 10 times for $\text{Test}_d$ (orange) and $\text{Train}_d$ (blue) reference objects with the following parameters: }
\begin{tabular}{c|c|c|c|c}
    prob & s & n &  homological degree & T \\ \hline
   uniform  &  $30$ & $2000$ & $0$ (left), $1$ (middle), $2$ (right) & $\infty$ 
\end{tabular}
\label{illus:adfhgsdfhfgh}
\end{figure}

As seen in Figure~\ref{illus:adfhgsdfhfgh}, the reference objects $\text{Test}_d$ and $\text{Train}_d$, for $d=2,7,8$  have noticeably different global stable ranks, indicating that there is some variation in the geometry between the training and test datasets. Since there is a probabilistic step in our pipeline, the whole process is repeated  10 times to demonstrate stability of the outcome.  

 A measure of geometric similarity between the reference objects can be obtained by considering distances between  the obtained stable ranks, for example by using the $L_1$ distance. We compute the average stable ranks corresponding to the training and test set respectively and present the distance between them, for each digit, in Table \ref{table:1}. To further investigate whether the difference corresponds to a dataset shift or is due to random factors we pool the training and the test set together and perform random partitions. This is done 10 times for each digit, average stable ranks are then computed and the distance between the training and test sets resulting from these random partitions is compared to the distances obtained for the original training and test split. The results indicate that the difference between training and test sets is not due to random factors alone. Perhaps it results from the way the dataset was originally partitioned (partitioned by writers, but several samples belong to each writer hence potentially introducing a bias).

 \begin{table}[h]
\begin{center}
\scalebox{0.85}{\begin{tabular}{c|c|c|c|c|c|c|c|c|c|c}
    Digit & 0 & 1 & 2 & 3 & 4 & 5 & 6 & 7 & 8 & 9 \\ \hline\hline
   H0 original split  &  356.3 & 365.3 & 1143.7 & 982.7 & 728.0 & 294.7 & 550.1 & 702.6 & 427.5 & 843.2 \\ 
   H0 random split & 122.3 & 93.8 & 125.7 & 134.3 & 41.6 & 119.8 & 42.2 & 101.0 & 132.6 & 57.1 \\ \hline
   H1 original split & 36.1 & 10.1 & 30.1 & 15.0 & 27.2 & 54.1 & 31.0 & 49.5 & 46.5 & 24.3 \\
   H1 random split & 7.7 & 3.4 & 6.7 & 3.6 & 1.7 & 13.5 & 2.4 & 2.3 & 9.4 & 4.8 \\ \hline
   H2 original split & 7.3 & 0.2 & 4.6 & 1.1 & 6.1 & 7.9 & 2.6 & 4.9 & 7.5 & 2.2 \\
   H2 random split & 1.6 & 0.6 & 0.6 & 0.8 & 0.4 & 3.3 & 0.5 & 0.3 & 0.8 & 0.8
\end{tabular}}
\end{center}
\caption{ Distances between stable ranks corresponding to training and test sets for different digits and homological degrees. The data is presented for stable ranks corresponding to the original training and test split and for stable ranks corresponding to random splits into training and test set. }\label{table:1}
\end{table}

\begin{figure}[h]
\begin{center}
\includegraphics[height=4cm]{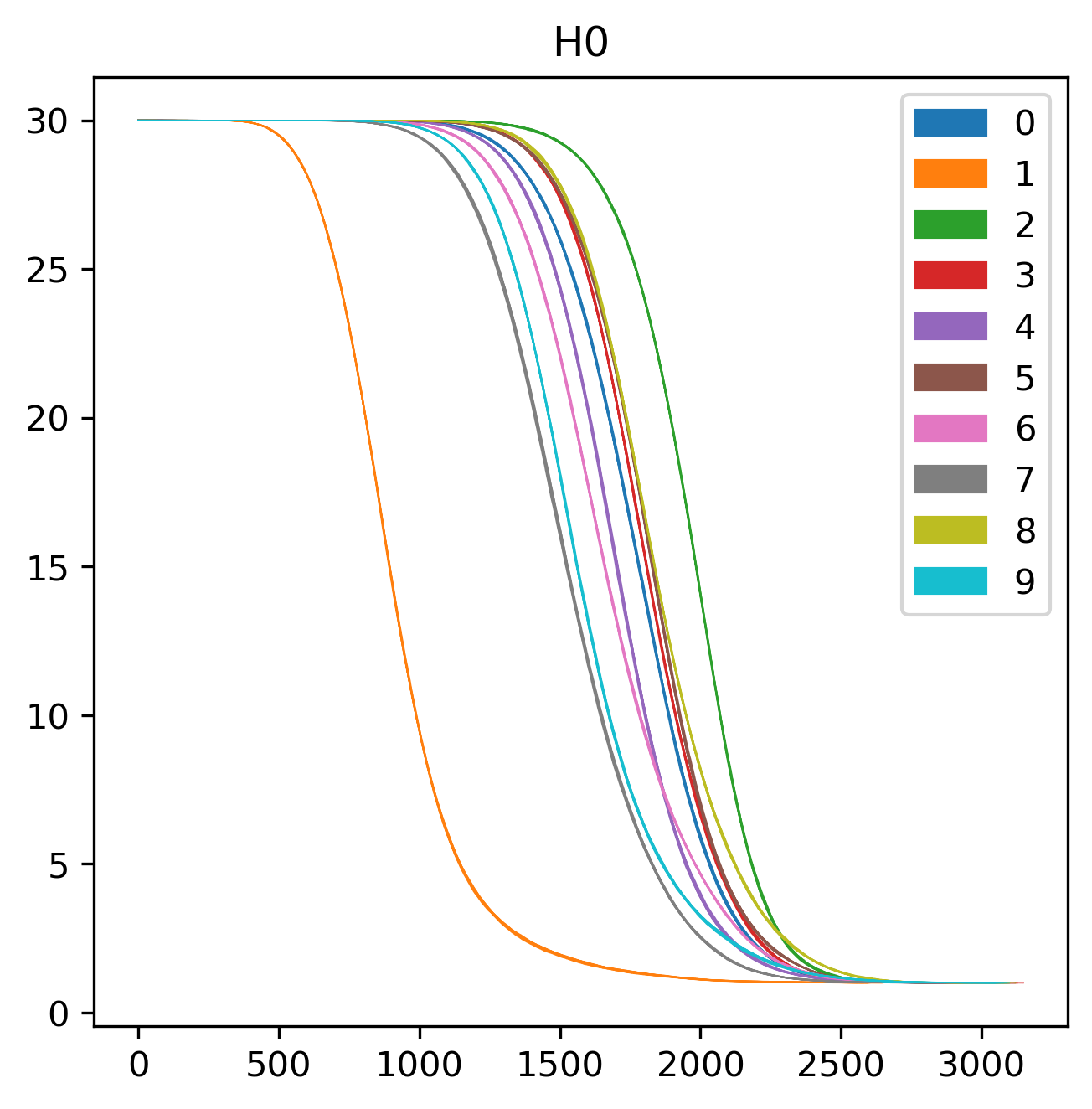}
\includegraphics[height=4cm]{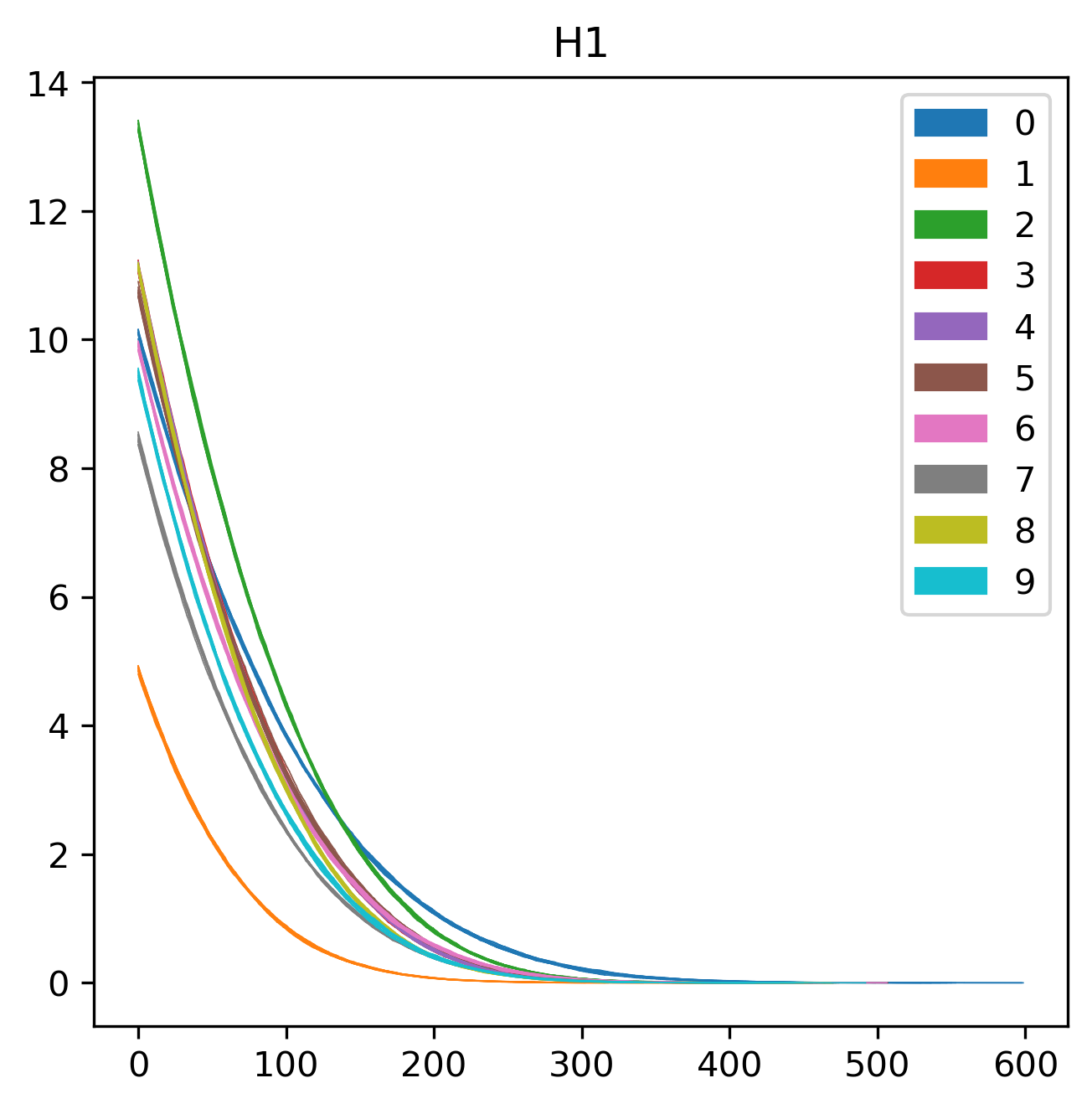}
\includegraphics[height=4cm]{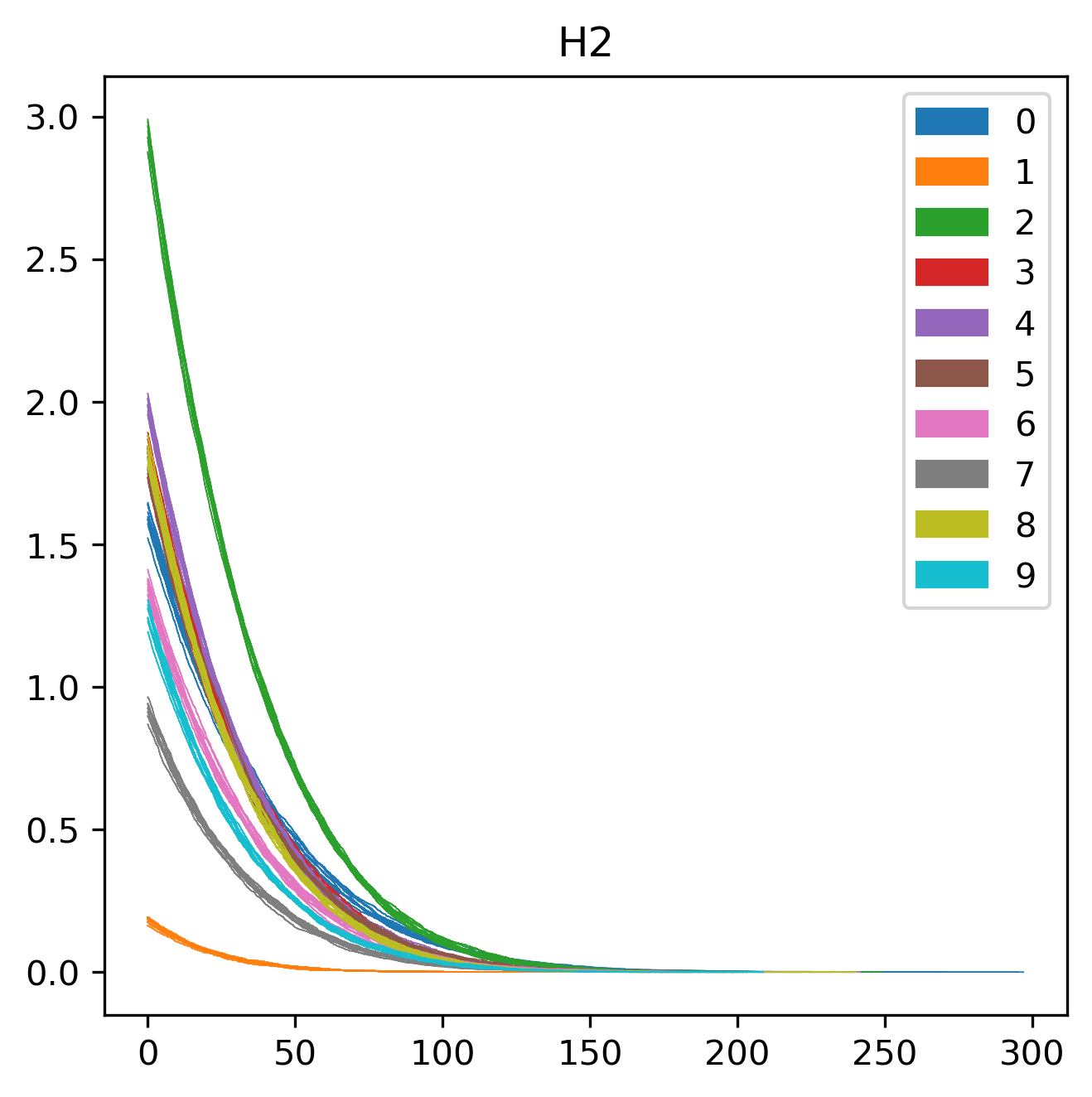}
\end{center}
\caption{The pipeline was repeated 10 times for every reference objects  $\text{Train}_d$, for all digits $d$, with the following parameters: }
\begin{tabular}{c|c|c|c|c}
    prob & s & n &  homological degree & T \\ \hline
   uniform  &  $30$ & $2000$ & $0$ (left), $1$ (middle), $2$ (right) & $\infty$ 
\end{tabular}
\label{illus:adfbsfbfgbhbfb}
\end{figure}

When we write a digit, we intuitively know which variations still enable communication.  We can think about the space $\text{Train}_d$ as a space encoding  such possible variations of $d$. A basic question is how dependent these spaces are on the digits and  whether these spaces, for different digits, have detectable global geometrical differences.
 Figure~\ref{illus:adfbsfbfgbhbfb} illustrates   some global stable ranks of these spaces.

 We note that in our experiments, the global stable ranks obtained by subsampling $s$-element subspaces were as (and sometimes more) distinctive of the digits as the stable ranks one can obtain from the computation of persistent homology on the whole reference object, without subsampling, a procedure that is heavier computationally (for instance, for homological degree 2 it was impossible to compute on a personal computer). 
 
\begin{figure}[h]
\centering
\includegraphics[height=4cm]{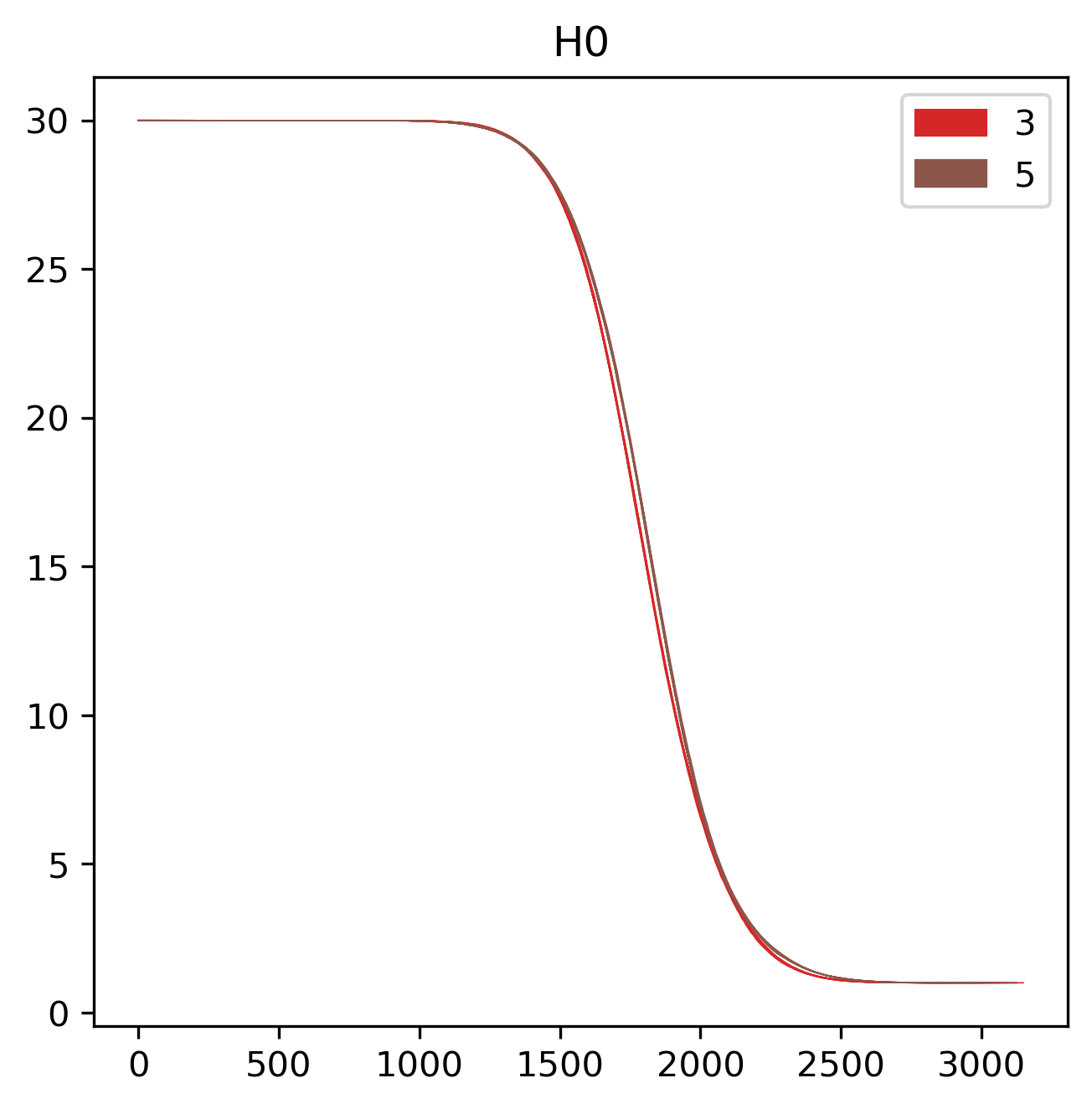}
\includegraphics[height=4cm]{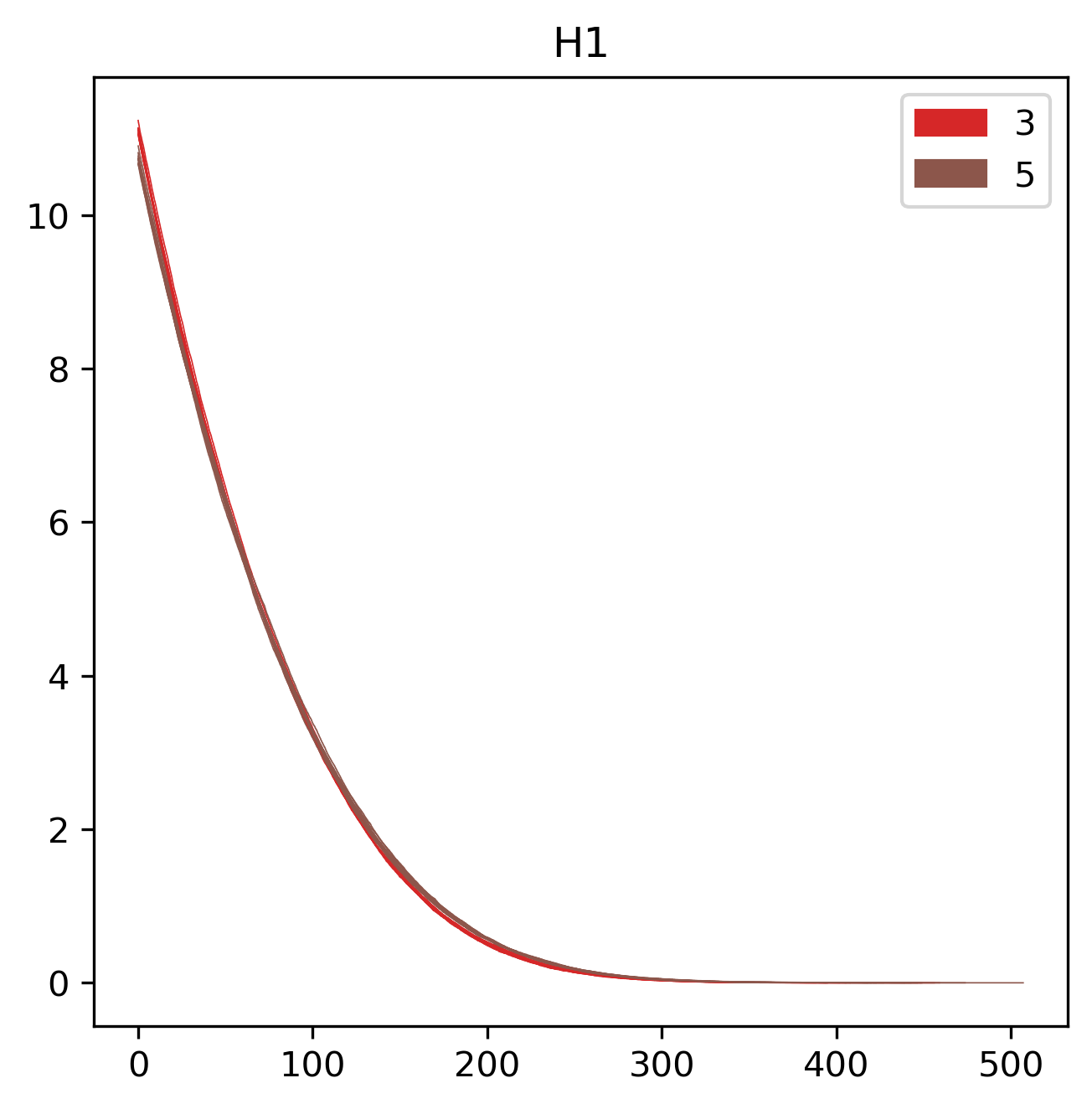}
\includegraphics[height=4cm]{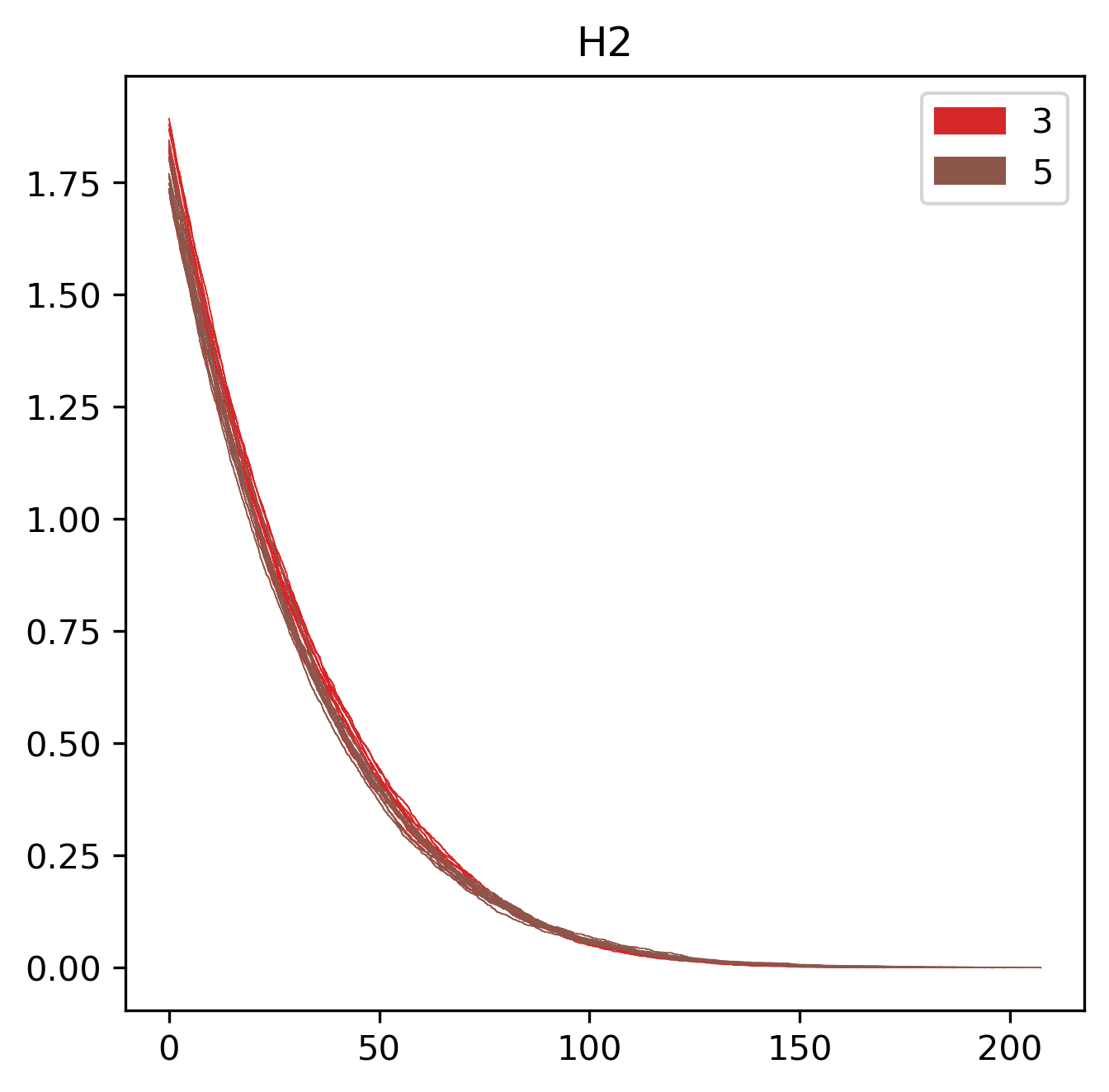}

\caption{The restriction of Figure~\ref{illus:adfbsfbfgbhbfb} to digits $3$ and $5$. }
\label{illus:asfdghshsfhfgh}
\end{figure}

In Section~\ref{sec:dfgsfghfg}, we discuss a strategy of how to use  the geometry of  the spaces  $\text{Train}_{d}$, encoded through our pipeline, to  classify  handwritten digits. Presented  examples suggest that the further apart the geometrical properties of the spaces  $\text{Train}_{d_1}$ and $\text{Train}_{d_2}$
are, the easier it is to distinguish between handwritten digits labeled  by $d_1$ and $d_2$. This indicates that we need to look for ways of amplifying geometrical differences, if there are any, between the spaces
$\text{Train}_{d}$ for various $d$.
For instance, consider the digits $3$ and $5$ and restrict Figure~\ref{illus:adfbsfbfgbhbfb}  to just these digits, illustrated in Figure~\ref{illus:asfdghshsfhfgh}.

\begin{figure}[h]
\begin{center}
\includegraphics[height=5cm]{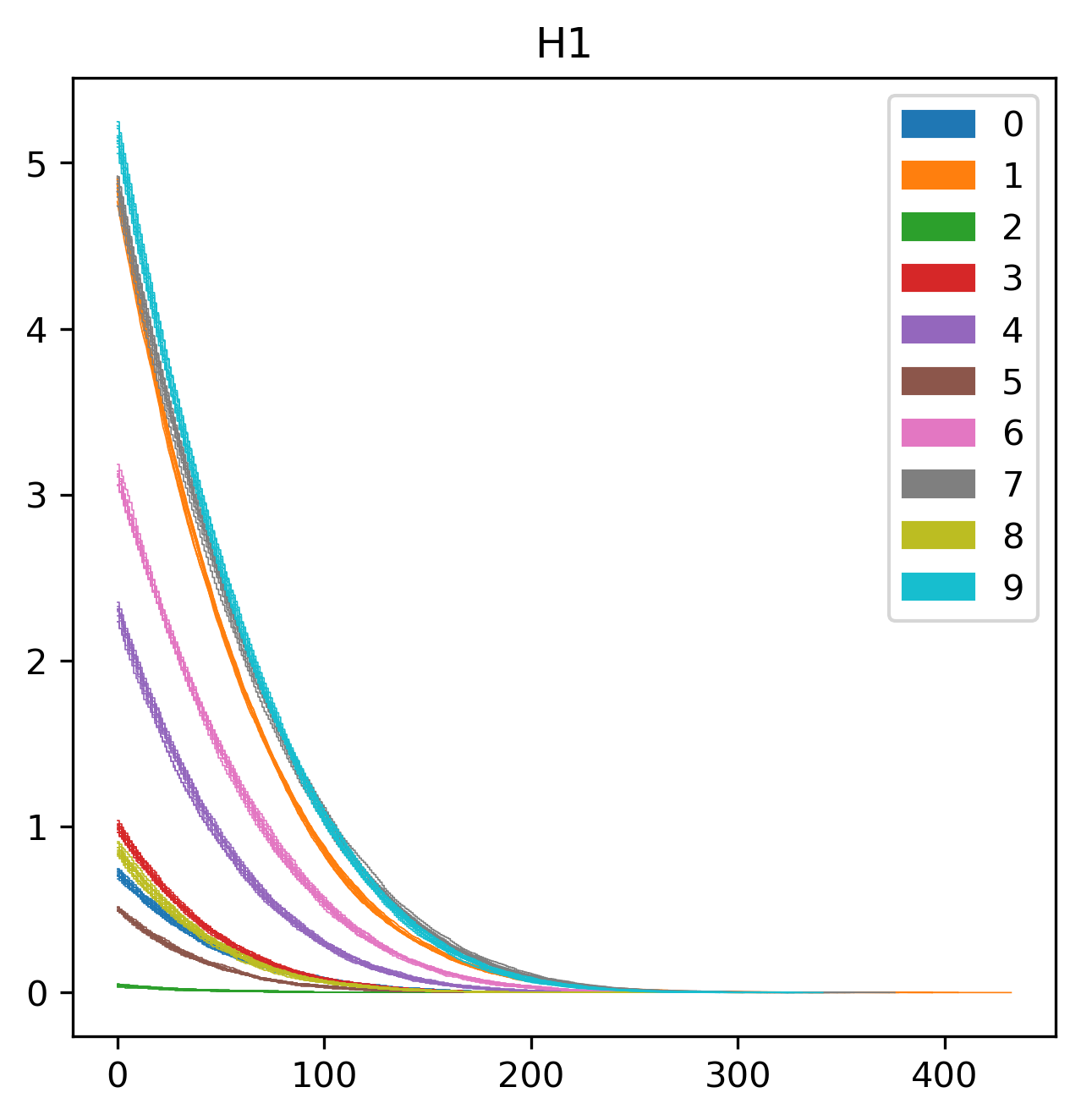}
\includegraphics[height=5cm]{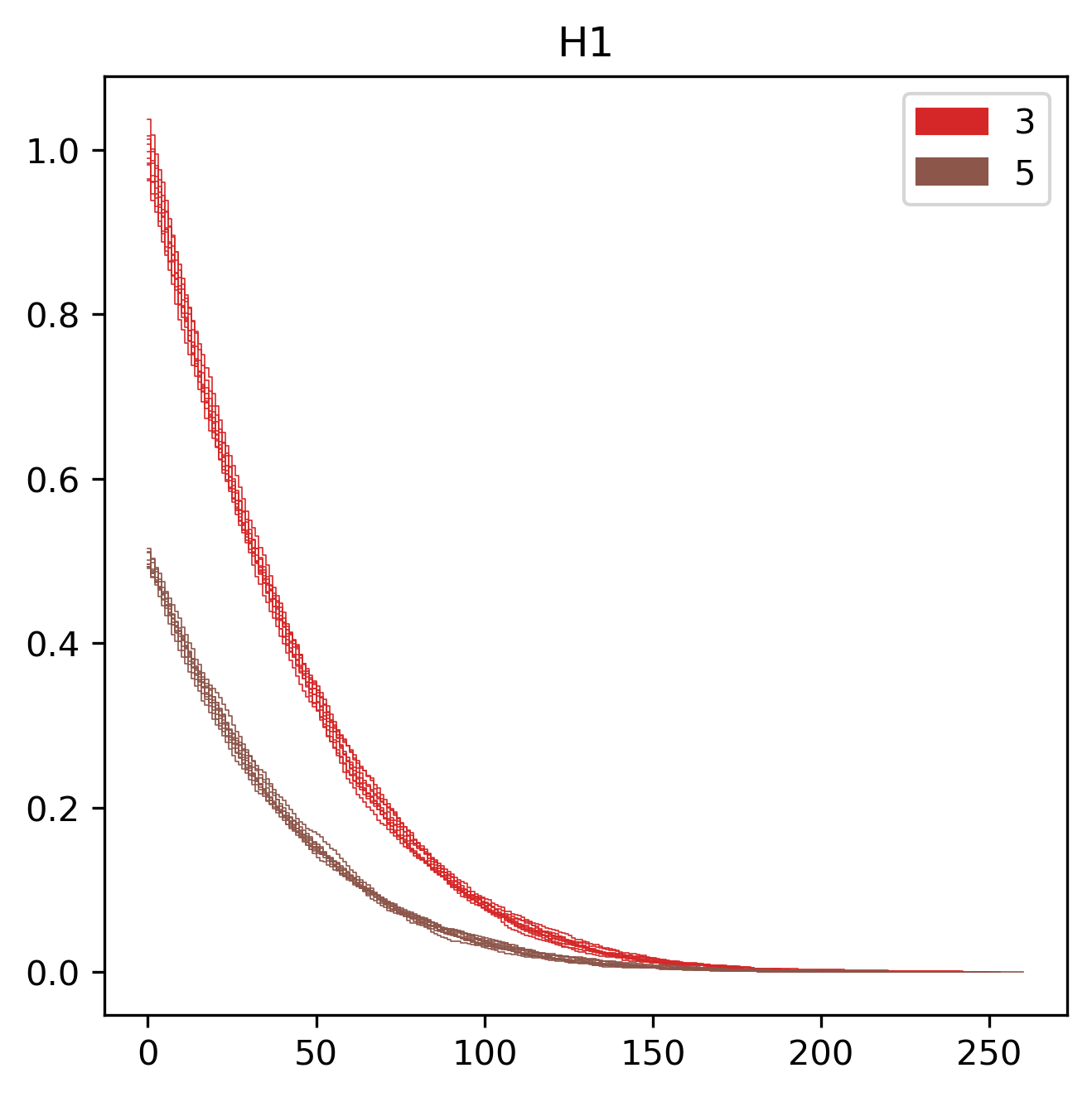}
\end{center}
\caption{The pipeline was repeated 10 times for every reference object  $\text{Train}_d$, for all digits $d$ (left), and for digits 3 and 5 (right), with the following parameters:}
\begin{tabular}{c|c|c|c|c}
    prob & s & n &  homological degree & T \\ \hline
   uniform  &  $30$ & $2000$ & $1$ & $1800$ 
\end{tabular}
\label{illus:dvgdfgbfdbhsdfgbh}
\end{figure}

Figure~\ref{illus:asfdghshsfhfgh} shows that our pipeline, with the
parameters used to produce this figure, does not see  geometrical differences between 
the spaces  $\text{Train}_{3}$ and $\text{Train}_{5}$. Are these spaces  then geometrically different and if so how can we encode differences between them? Let us change the truncation parameter $T$ to $1800$. The effect is shown in Figure~\ref{illus:dvgdfgbfdbhsdfgbh}, illustrating the fact that varying the parameters, e.g. sample size or the parameters used to construct the stable ranks, can lead to stable distinctive descriptors.

\section{Relative stable ranks in the plane}
\label{sec:asgwdfhfg}
The results of the pipeline described in Section~\ref{sec:pipeline}, when the outcome of step A is given by the  probability function determined by a point $p$, are called {\em relative stable ranks} of  the reference object. We think about  relative stable ranks  as encoding geometrical information about the position of the point $p$ in the ambient space $\mathbb{R}^r$ in relation to the reference object. In this section we illustrate how  relative stable ranks can be used to describe simple geometrical aspects of points in $\mathbb{R}^r$.   Our initial data $X$ consist of 200 random points on the plane (consisting of both orange and blue points in  Figure~\ref{illus:sdfbgsfghgvg})  whose positions we would like to geometrically describe. 

 \subsection*{Example 1}
 
\begin{figure}[h]
\begin{center}
\includegraphics[height=5.5cm]{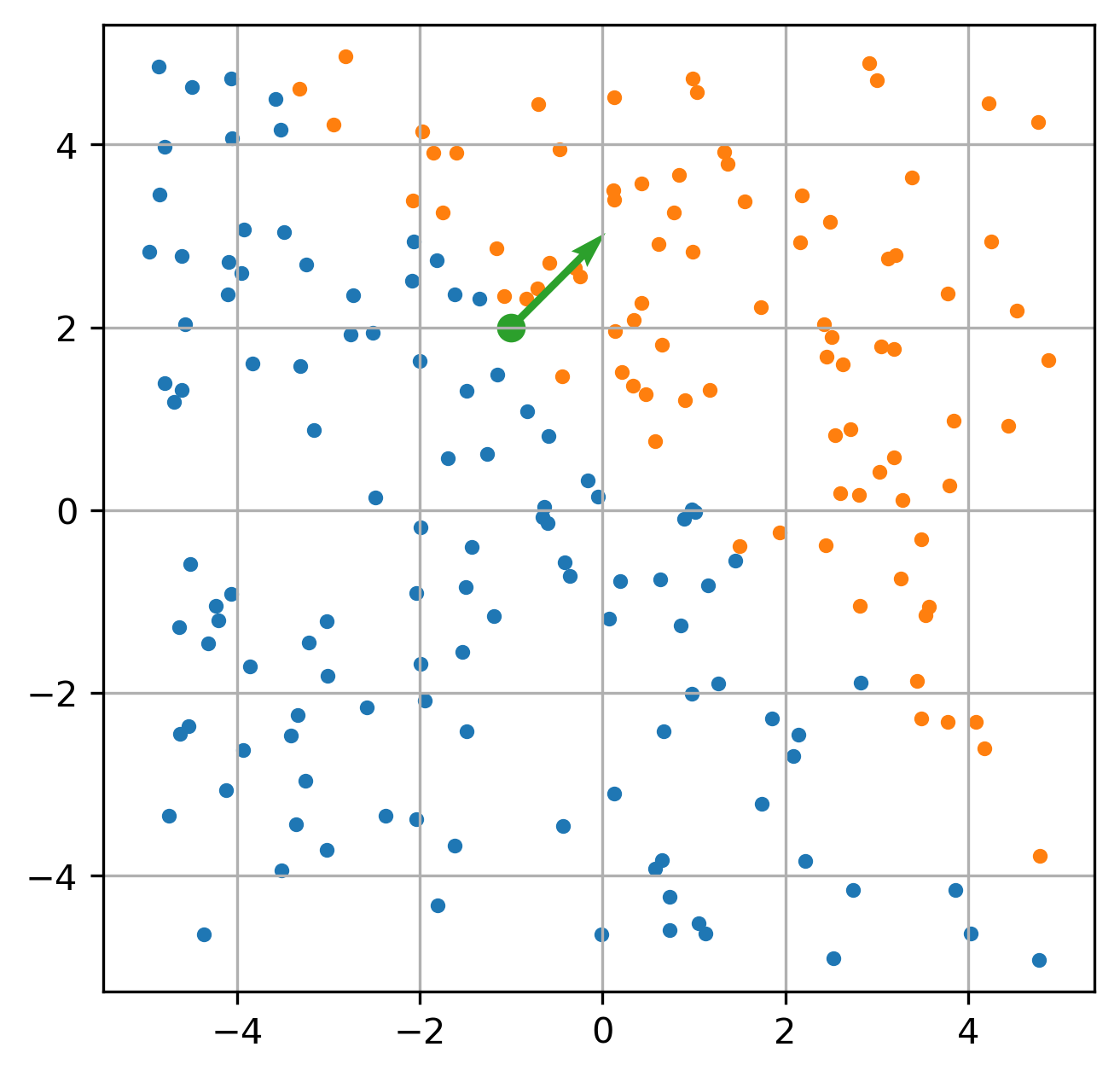}
\includegraphics[height=5.5cm]{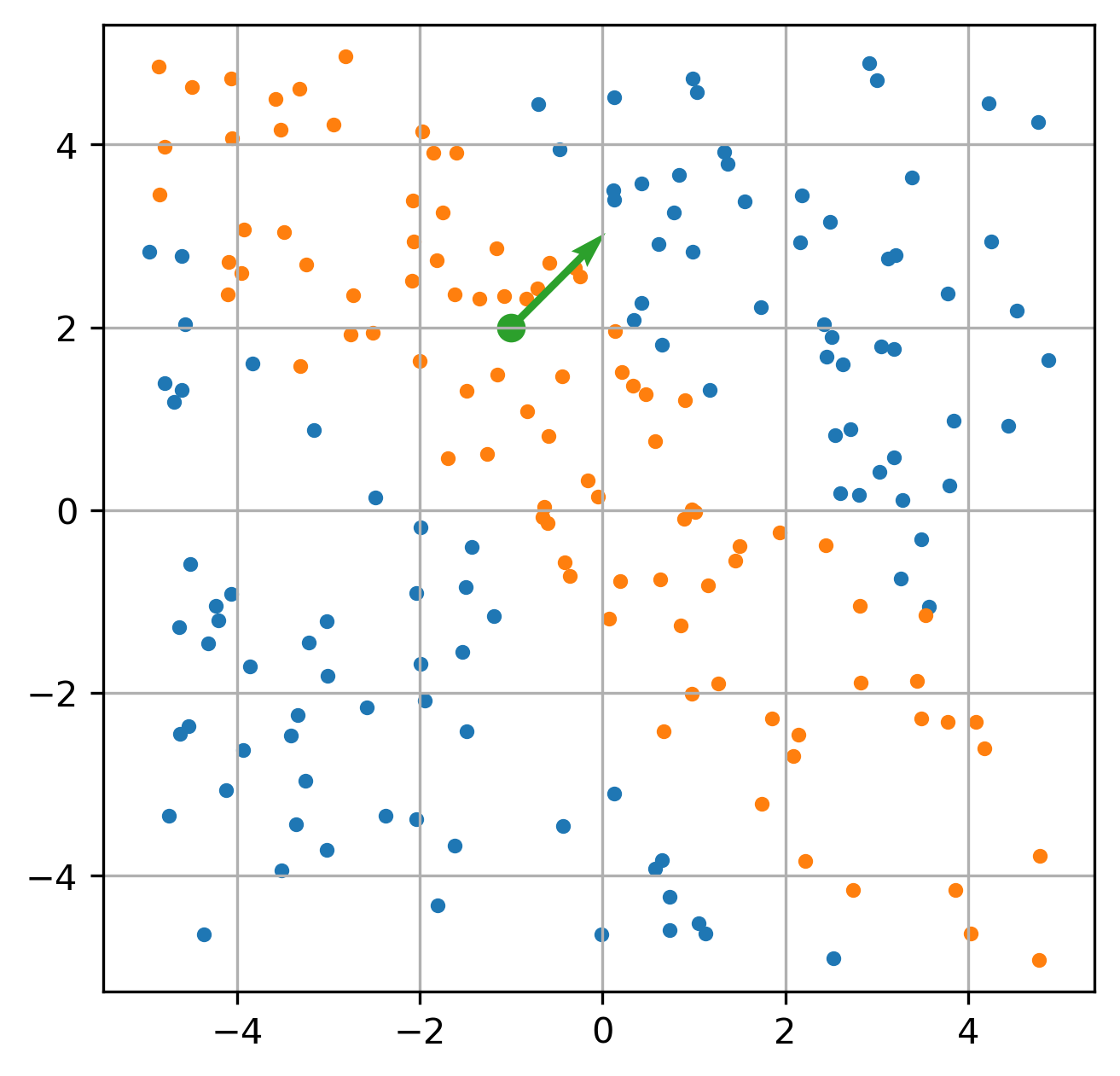}
\end{center}
\caption{ 200 random points colored according to whether the corresponding stable rank has constant value $1$ (orange) or $0$ (blue). The stable ranks were obtained with reference object containing only  point $[-2,1]$ and the vector field given by the vector $[1,1]$ (Example 1 Section \ref{sec:asgwdfhfg}).
}\label{illus:sdfbgsfghgvg}
\end{figure}

\begin{itemize}
    \item {\em Reference object}: a single point whose coordinates are $[-1,2]$.
    \item {\em Point}: any point $p$ in  $X$. 
    \item {\em Vector field}: given by the vector $[1,1]$.
    \item {\em Distribution}: we consider two distributions:
    \[\mathcal{D}_1(x) := \begin{cases} 1 &\text{ if }x\geq 0\\
    0& \text{ if } x<0
    \end{cases}\ \ \ \ \ 
    \mathcal{D}_2(x) := \begin{cases} 1 &\text{ if }-2\leq x\leq 1\\
    0& \text{ otherwise }
    \end{cases}\]
    \item {\em The other parameters}: $s=1$, $n=1$, homological degree is $0$, and $T=\infty$.
\end{itemize}
Since the reference object consists of just one point,
other  homological degrees than $0$ are irrelevant.  
The outcome of our pipeline in this case, for every point $p$ in $X$,
 is a constant function $0$ or $1$. In this way the initial dataset $X$  is partitioned into two clusters: points leading to the  stable rank $0$ and points  leading to the  stable rank  $1$. The two illustrations in Figure~\ref{illus:sdfbgsfghgvg}, which correspond to
 the  two distributions $\mathcal{D}_1$ and $\mathcal{D}_2$, 
 show such  partitions of $X$.
 We see that our pipeline  can for example distinguish  between points lying on different sides of a hyperplane, an interesting piece of geometrical information.

\subsection*{Example 2}

\begin{figure}[h]
\begin{center}
\includegraphics[height=5.5cm]{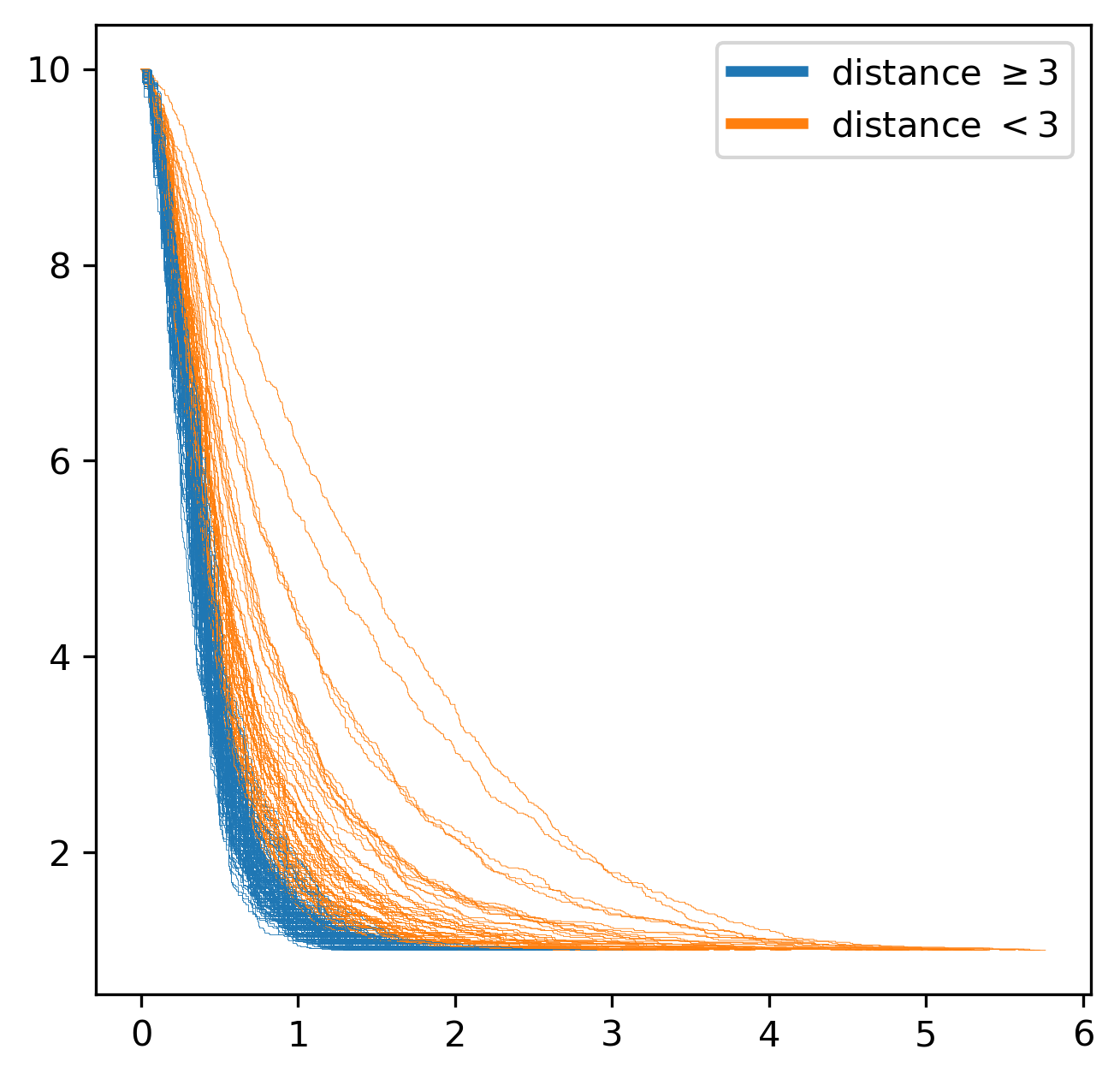}
\includegraphics[height=5.5cm]{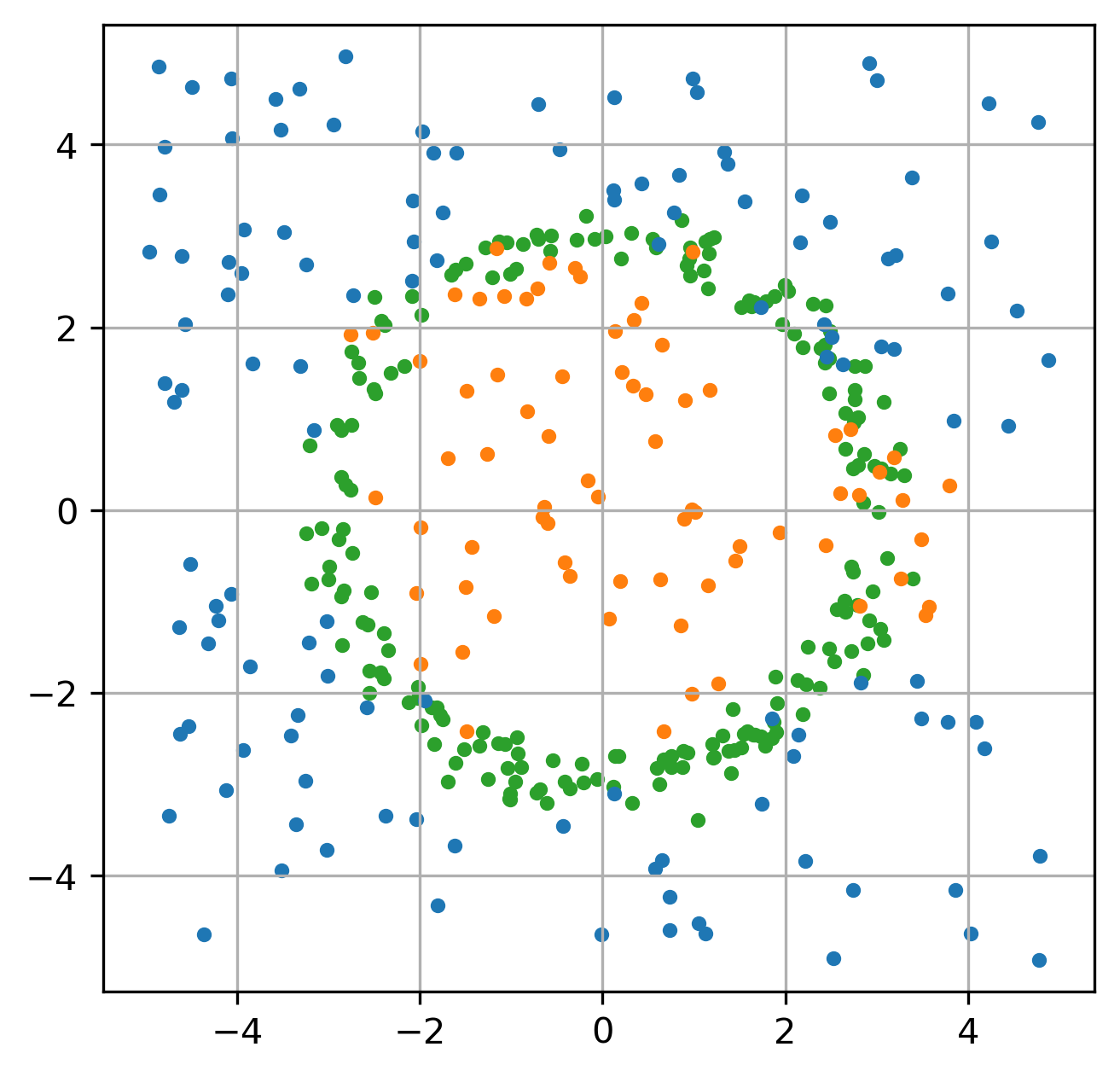}
\end{center}
\caption{ \textbf{Left}: Stable ranks corresponding to the random points in the plane, colored according to their distance to the origin. \textbf{Right}: Reference object (green) and random points colored according to whether the corresponding stable rank at  $0.87$ has value bigger (orange) or lower (blue) than $1.87$ (Example 2 Section \ref{sec:asgwdfhfg}). }\label{illus:hrtjkryukj}
\end{figure}

\begin{itemize}
    \item {\em Reference object}: a noisy circle (of radius $3$) represented  by green dots in Figure~\ref{illus:hrtjkryukj}.  
    \item {\em Point}: any point $p$ in  $X$. 
    \item {\em Filter}: assigns to an element in the reference object its distance to $p$.
   \item {\em Distribution}: Gaussian centered at $0$ with standard deviation $1$.
    \item {\em The other parameters}: $s=10$, $n=100$, homological degree is $0$, and $T=\infty$.
\end{itemize}
In Figure~\ref{illus:hrtjkryukj} on the left, obtained stable ranks for all points in $X$ are plotted. Those stable ranks corresponding to points  whose distance to the origin is less than 3 are orange and the others are blue. In the illustration on the right a point is orange if corresponding stable rank at  $0.87$ has value  bigger than $1.87$. The other points of $X$ are blue. Green dots represent the reference object. In this case we see that our pipeline can be used to decide if a point is inside or outside  a circle, again an interesting geometrical property. In this case a simple rule allowed us to discriminate between points inside and outside the circle, based on their stable ranks. In the next section we will see that such classification rules can also be learned from the data.

Deciding if a points is inside or outside a circle can be obtained by our pipeline with another set of  parameters:   

\subsection*{Example 3}

\begin{figure}[h]
\begin{center}
\includegraphics[height=5.5cm]{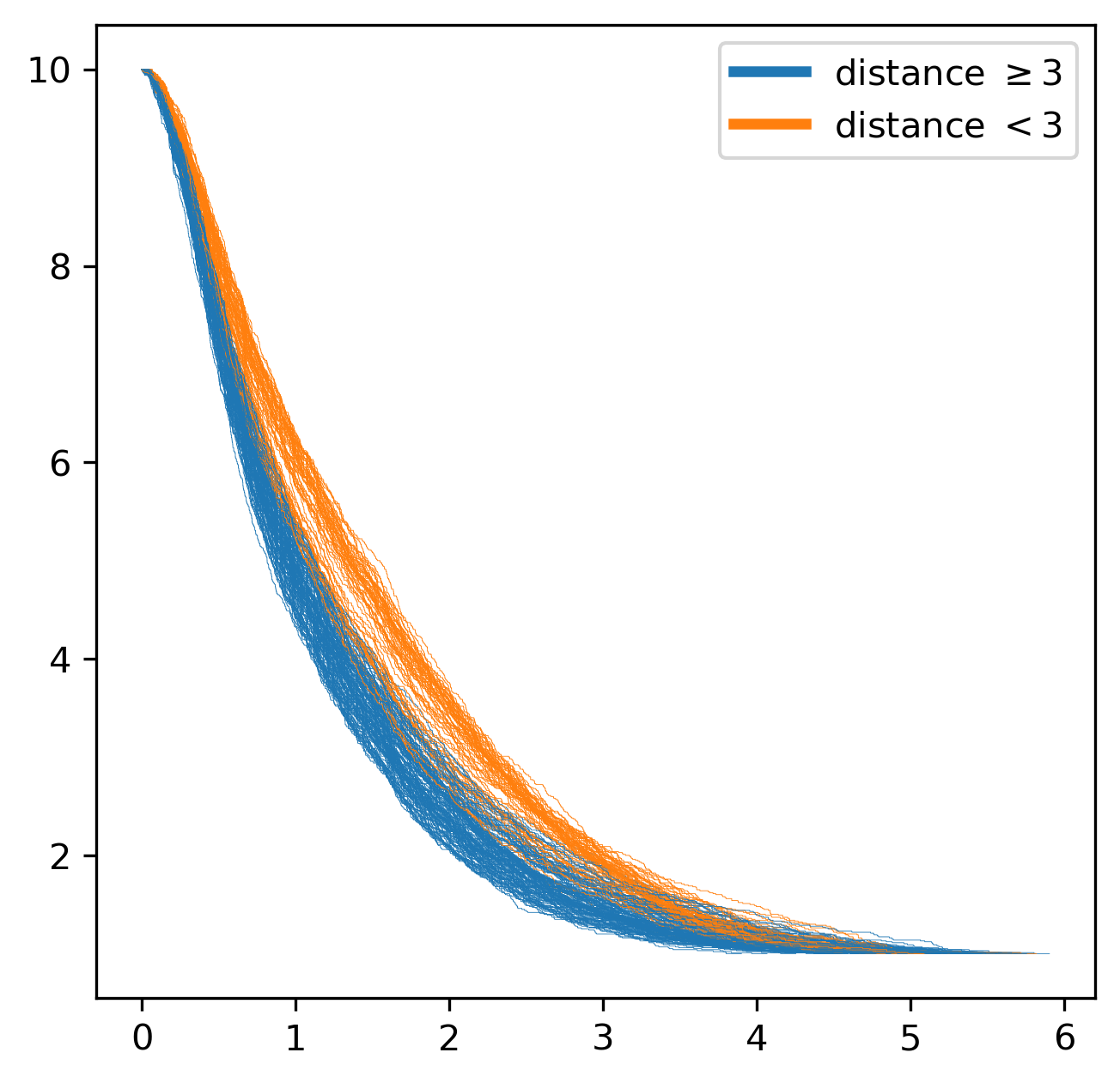}
\includegraphics[height=5.5cm]{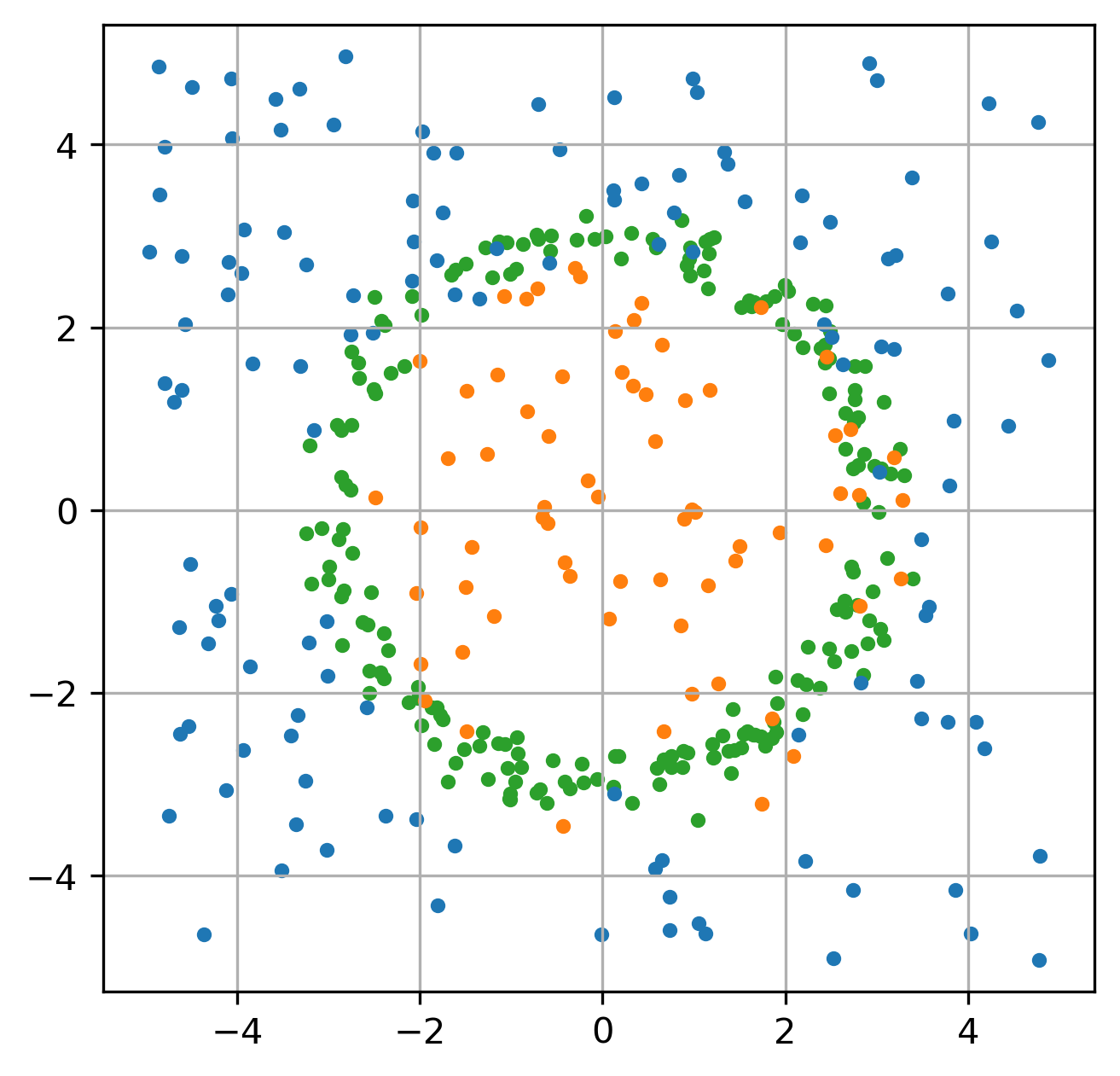}
\end{center}
\caption{ \textbf{Left}: Stable ranks corresponding to the random points in the plane, colored according to their distance to the origin. \textbf{Right}: Reference object (green) and random points colored according to whether the corresponding stable rank at  $1.5$ has value bigger (orange) or lower (blue) than $3.97$ (Example 3 Section \ref{sec:asgwdfhfg}). }\label{illus:qwerrtrtyujt}
\end{figure}

\begin{itemize}
    \item {\em Reference object}: the same noisy circle as in Example 2.  
    \item {\em Point}: any point $p$ in  $X$. 
    \item {\em Vector field}: assigns to an element in the reference object the
    vector from that element to the center of mass of the reference object (which in this case is close to the origin).
   \item {\em Distribution}: $\mathcal{D}(x) := \begin{cases} 1 &\text{ if }  x\geq 1\\
    0& \text{ if } x<1
    \end{cases}$.
    \item {\em The other parameters}: $s=10$, $n=100$, homological degree is $0$, and $T=\infty$.
\end{itemize}
In Figure~\ref{illus:qwerrtrtyujt} on the left, obtained stable ranks for all points in $X$ are plotted. As in Example 2, those stable ranks corresponding to points, whose distance to the origin is less than 3, are orange and the other are blue. In the  illustration on the right  a point is orange if the corresponding stable rank  at $1.5$ has value  bigger than $3.9$. The other points of $X$ are blue, and the green dots represent the reference object. We see again that our pipeline can be used to decide if a point is inside or outside  a circle.


\section{Relative stable ranks on MNIST}
\label{sec:dfgsfghfg}

In Section \ref{sec:global} we provided experiments on MNIST for global stable ranks. In this section we now shift the focus to relative stable ranks. In the following experiments, subsets of the training sets corresponding to one or several digits will be used as reference objects. We will illustrate that these reference objects, when sampled from the perspective of different types of points in $\mathbb{R}^{784}$, such as points corresponding to digits in the test set, have interesting geometries.

Following the steps defined in the pipeline, for a point under consideration $p$ and for elements in the reference object $x$ in $\mathcal{R}$ we choose as \textit{filter function} $f_p(x) = ||p-x||_{2}$, i.e. the Euclidean distance between the point under consideration and the elements of the reference object. As \textit{distribution} we choose a Gaussian, whose parameters $\mu_p, \sigma_p$ are chosen in order to concentrate the probability mass on elements of the reference object close to $p$, yet ensuring the probability mass is distributed on sufficiently many elements for the samples to be diverse enough. We consider the set $\text{Dist}_p$ of all distances between $p$ and points in $\mathcal{R}$, i.e. $f_p(x)$ for all $x \in \mathcal{R}$. We select $\mu_p$ to be the k:th percentile of $\text{Dist}_p$, where k typically is a low number. We then choose $\sigma_p$ in relation to the \textit{sample size} parameter such that \textit{sample size} $\times$ \textit{amplification} elements of $\text{Dist}_p$ lie within one standard deviation, where \textit{amplification} is also a fixed parameter.

The distances between the point under consideration and the points in the reference object are only used indirectly, in order to construct a probability distribution on the reference object, the geometry of which is our interest. Moreover, we note that two filter functions $f_p$ and $C f_p$, $C \in (0, \infty)$ yield the same probability distribution on $\mathcal{R}$, that is we are only interested in the relative distances. We believe that this makes our approach fundamentally different from distance-based methods such as KNN.

\subsection{Illustration of the pipeline and first example}

We start with a basic example to illustrate the pipeline. We take as our reference object 
the set $\text{Train}_1$,
of all samples from the MNIST training set corresponding to the digit 1. 
Next we select two points from the ambient space, $\mathbb{R}^{784}$: the origin of that space and the center of mass of the reference object.

\begin{figure}[h]
\begin{center}
\includegraphics[width=12.5cm]{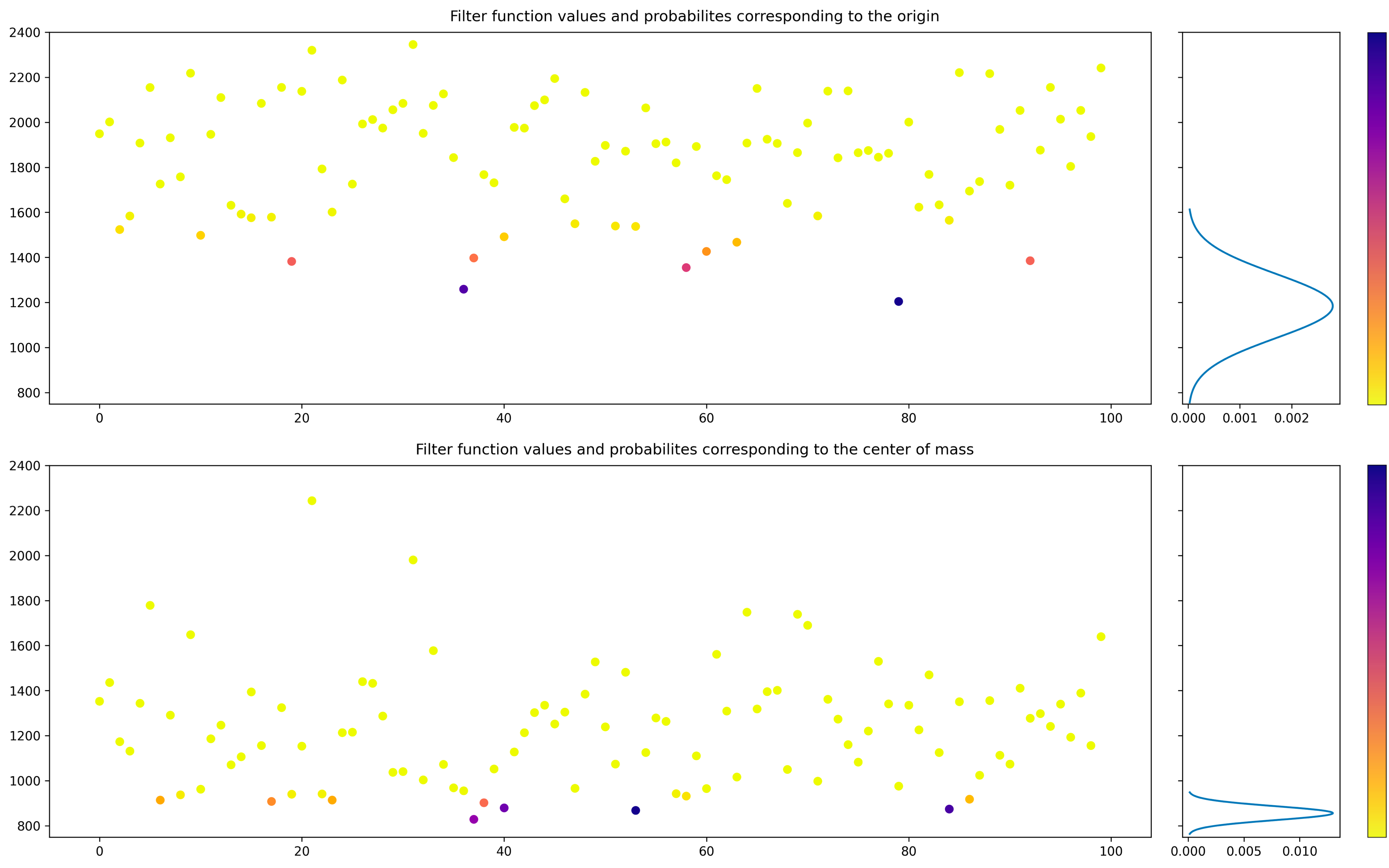}
\end{center}
\caption{Filter function values corresponding to the origin (top) and center of mass (bottom) for 100 elements of the reference object. Dots are colored according to their probability. To the right, the Gaussian corresponding to the point.}\label{illus:0}
\end{figure}

In Figure \ref{illus:0}, for the clarity of the plot we select 100 elements out of the 6742 composing the reference object. The y-axis corresponds to the value of the filter function for those elements, that is the distance to the origin (top plot), and the distance to the center of mass (bottom plot). Based on the values of the filter function, a Gaussian is chosen for each of the two points and shown to the right of the subplots (we use as parameters \textit{k:th percentile}$=1$, \textit{amplification}$=5$). Next, as described in Section \ref{sec:pipeline:a2}, a probability distribution on the reference object is computed for each point, by evaluating the values of the filter function under the Gaussian and normalizing. In the plot, the probabilities of the elements of the reference object are indicated by the color of the dot. By comparing the subplots we see that different probabilities are attached to the elements of the reference object.

\begin{figure}[h]
\begin{center}
\includegraphics[width=8cm]{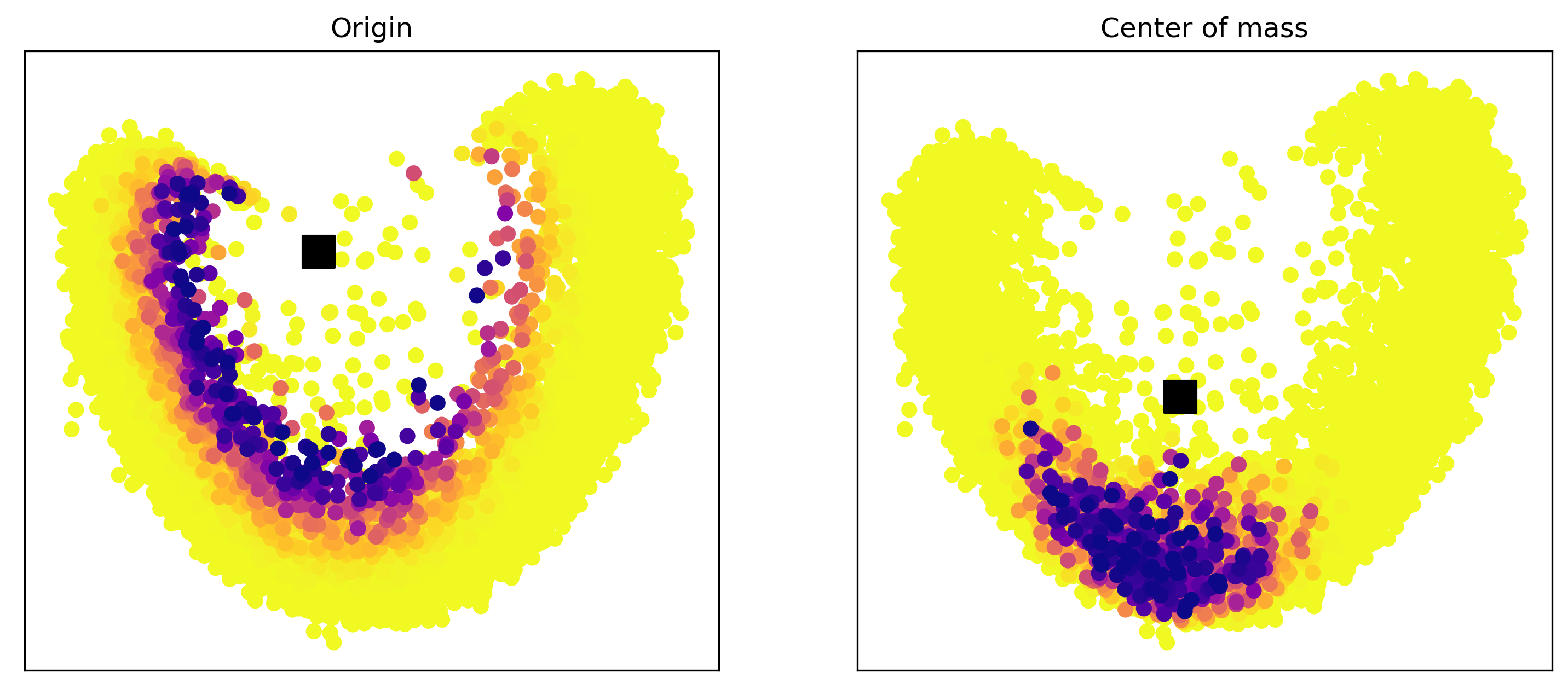}
\end{center}
\caption{Projection on the two first principal components of the reference object. The origin (left) and the center of mass (right) are represented by black squares. Other dots are colored according to their probability. }\label{illus:1}
\end{figure}

We further illustrate this idea in Figure \ref{illus:1}. The two first principal components of the reference object are computed. We then project the reference object together with the origin and center of mass on the principal components. The origin (left plot) and the center of mass (right plot) are represented by black squares, and the dots representing elements of the reference object are colored according to their probability in the same way as in the previous plot.

\begin{figure}[h]
\begin{center}
\includegraphics[width=12.5cm]{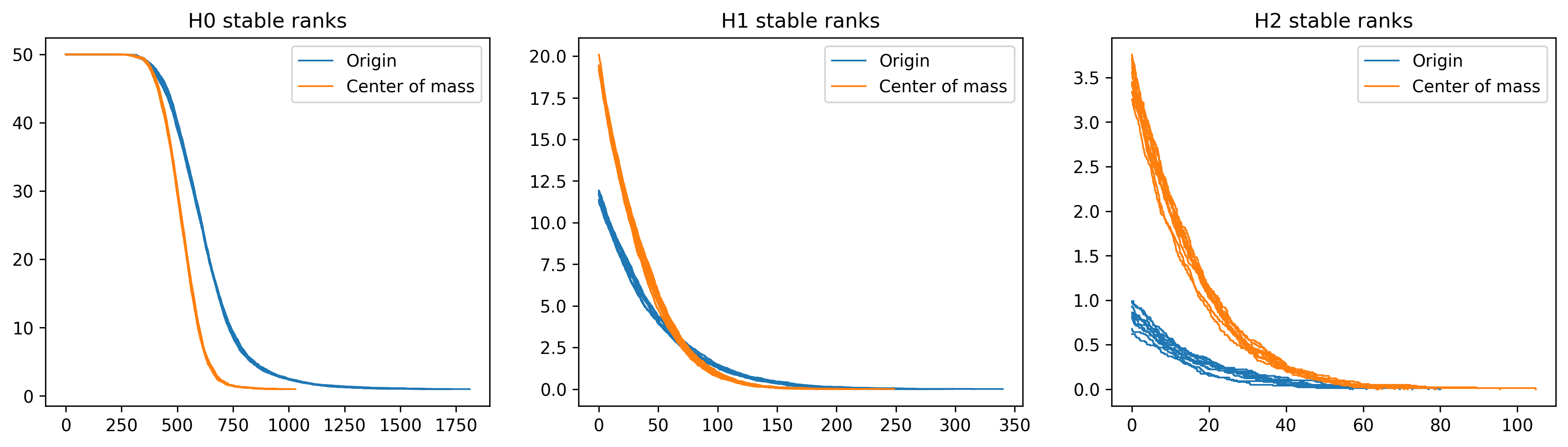}
\end{center}
\caption{Stable ranks corresponding to the origin and the center of mass for different homological degrees.}\label{illus:2}
\end{figure}

Having illustrated that different points lead to different probability distributions, we now subsample the reference object according to these probability distributions. For each such sample a distance space is constructed (with Euclidean distance). Next, as described in Section \ref{sec:pipeline:b2}, these distance spaces are converted into persistence modules corresponding to each homological degree, and then to stable ranks (we use \textit{sample size}$=50$, \textit{number of instances}$=100$). The resulting average stable ranks, presented in Figure \ref{illus:2}, demonstrate that the geometrical signatures corresponding to the origin and the center of mass are distinct. We plot 10 stable ranks for each point and homological degree, obtained by repeating the whole procedure each time.

\begin{figure}[h]
\begin{center}
\includegraphics[width=14cm]{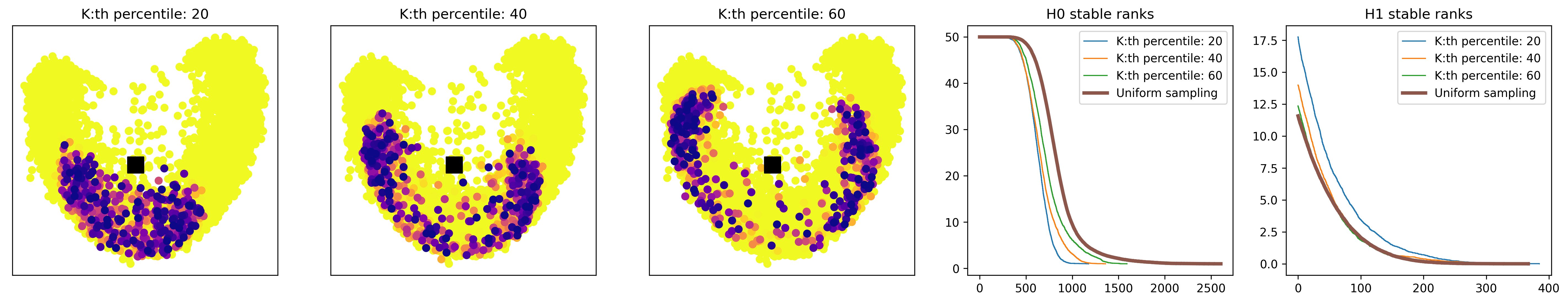}
\end{center}
\caption{\textbf{Plot 1, 2, 3}: Projection on the two first principal components of the reference object. The center of mass is represented by a black square. Other dots are colored according to their probability when varying the \textit{k:th percentile} parameter. \textbf{Plot 4, 5}: Stable ranks corresponding to the center of mass for the different parameters and stable rank corresponding to uniform subsampling, for homological degrees H0 and H1. }\label{illus:varying_perc}
\end{figure}

\begin{figure}[h]
\begin{center}
\includegraphics[width=14cm]{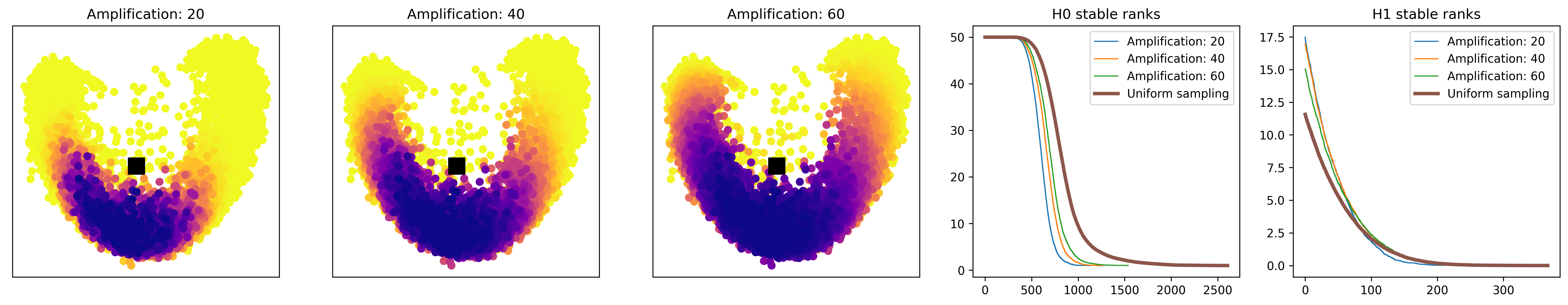}
\end{center}
\caption{\textbf{Plot 1, 2, 3}: Projection on the two first principal components of the reference object. The center of mass is represented by a black square. Other dots are colored according to their probability when varying the \textit{amplification} parameter. \textbf{Plot 4, 5}: Stable ranks corresponding to the center of mass for the different parameters and stable rank corresponding to uniform subsampling, for homological degrees H0 and H1.}\label{illus:varying_amplif}
\end{figure}

Another way to illustrate how the geometry changes as we subsample the reference object in different fashions is to take the perspective of one point only – here we choose the center of mass – but to vary the parameters defined in the previous section. In Figure \ref{illus:varying_perc} we show how the probability distributions on the reference object (illustrated on the PCA projection of the elements) and the resulting stable ranks change when we increase the \textit{k:th percentile} of the vector of distances on which our Gaussian is centered (all other parameters are held fixed), meaning that we sample elements at increasing distances to our point. In Figure \ref{illus:varying_amplif} we instead show the effect of increasing the \textit{amplification} parameter, which means that less probability mass will be concentrated on elements whose distance is close to the mean of the Gaussian. For both parameters, as their values increase, the stable ranks become closer to the global descriptor of the reference object described in Section \ref{sec:global}, i.e. to the stable rank obtained by uniform subsampling of the reference object, indicating that the geometry is less and less informative.

\subsection{Inside and outside}

Instead of choosing the center of mass of the whole reference object, we now perform k-means clustering (k=10) on the reference object and select the center of mass of each cluster. We also sample 10 points randomly from the ambient space (the subset of $\mathbb{R}^{784}$ corresponding to allowed pixel values). We can then apply the procedure described in the previous section to obtain 10 stable ranks for the points corresponding to the centers of mass and 10 stable ranks corresponding to the random points, for each homology degree.

\begin{figure}[h]
\begin{center}
\includegraphics[width=12.5cm]{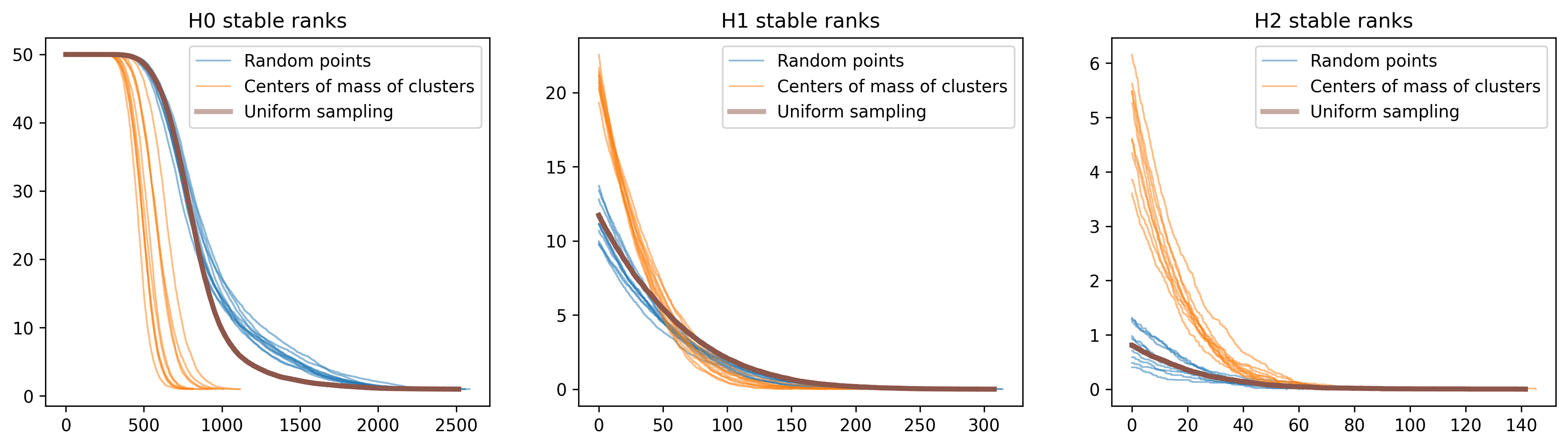}
\end{center}
\caption{Stable ranks corresponding to random points, to the centers of mass of clusters and to uniform sampling, for different homological degrees. }\label{illus:3}
\end{figure}

These stable ranks are displayed in Figure \ref{illus:3} together with the average stable rank corresponding to a uniform subsampling of the reference object. Our aim is to illustrate that stable ranks resulting from sampling from "inside" the reference object, e.g. for centers of mass, are distinct from stable ranks obtained by sampling from the "outside", e.g. from random points in the ambient space or from the origin (in the previous example). The latter are in turn more similar to the stable rank obtained by uniform subsampling.

\subsection{Distinguishing out-of-sample points from two subsets of the reference object}\label{expmnist:unionoftrainingsets}

In the previous section, we saw a clear distinction between stable ranks obtained by sampling from the "inside" and from the "outside" of the reference object. But the stable ranks corresponding to different centers of mass also displayed some variability, indicating a difference in the geometry. We aim to explore this idea further in the following setting: we now take as our reference object the union $\text{Train}_1\cup \text{Train}_7$ of all samples from the MNIST training set corresponding to digit $1$ or digit $7$. We note that we still have only one reference object and the labels (indicating whether an element of the reference object corresponds to a $1$ or a $7$) are not used. Instead of considering random points or centers of mass as in the previous section, we now consider 10 points randomly chosen from $\text{Test}_1$ and 
$10$ points randomly chosen from $\text{Test}_7$  
and repeat the same procedure to compute the stable ranks representing these points. In our pipeline we use the following parameters:  {\em sample size}$=30$, {\em amplification}$=2$, for the {\em  homological degree} $0$,  
the {\em truncation} parameter is set to be $\infty$,
and for the {\em  homological degree} $1$,   
the {\em truncation} parameter is set to be $1200$.

\begin{figure}[h]
\begin{center}
\includegraphics[width=11cm]{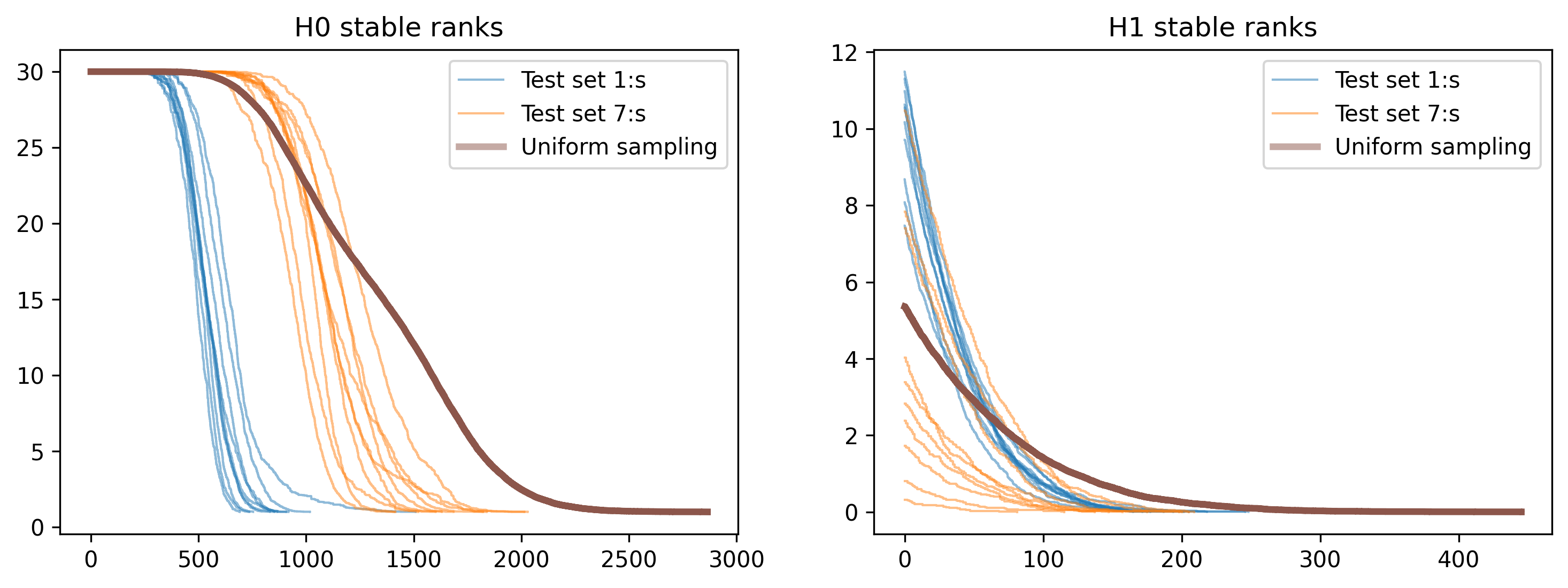}
\end{center}
\caption{Stable ranks corresponding to test set 1:s, test set 7:s and to uniform sampling, for different homological degrees. Training set 1:s and 7:s used as reference object.}\label{illus:4}
\end{figure}

We can see in Figure \ref{illus:4} that the stable ranks corresponding to test set digit 1 are distinct from those corresponding to test set digit 7, and they are both distinct from the stable ranks resulting from the uniform subsampling of the reference object. Hence, when we sample based on distances to test set digits 1:s or 7:s, we sample subsets of the reference object where the geometry is different, which allows us to discriminate between the points we sampled from.

To quantify the capacity to discrimate between digits based on their stable ranks, we train a Support-vector machine classifier on the 20 stable ranks, for each homological degree, using the kernel obtained by taking inner products between stable ranks in the $L_2$ function space \cite{10.3389/fams.2021.668046}. We can then evaluate the model on the remaining samples of digit 1 and 7 from the MNIST test set (samples that are neither part of the reference object nor part of the 20 samples used for the training). We obtain an accuracy of 96.9\% for H0 stable ranks and 94.5\% for H1 (average accuracy after repeating the procedure 10 times with different samples used for training). While we are not aiming at approaching state of the art accuracy levels we believe the results point to the fact that the geometry of a reference object, when chosen judiciously and in relation to a point, can be informative about characteristics of this point. We also note that we used a large unlabeled dataset (our reference object) but only a few (20) labeled samples, which is the setting of semi-supervised learning.

\subsection{Distinguishing out-of-sample points based on another reference object}\label{expmnist:otherrefobject}

\begin{figure}[h]
\begin{center}
\includegraphics[width=11cm]{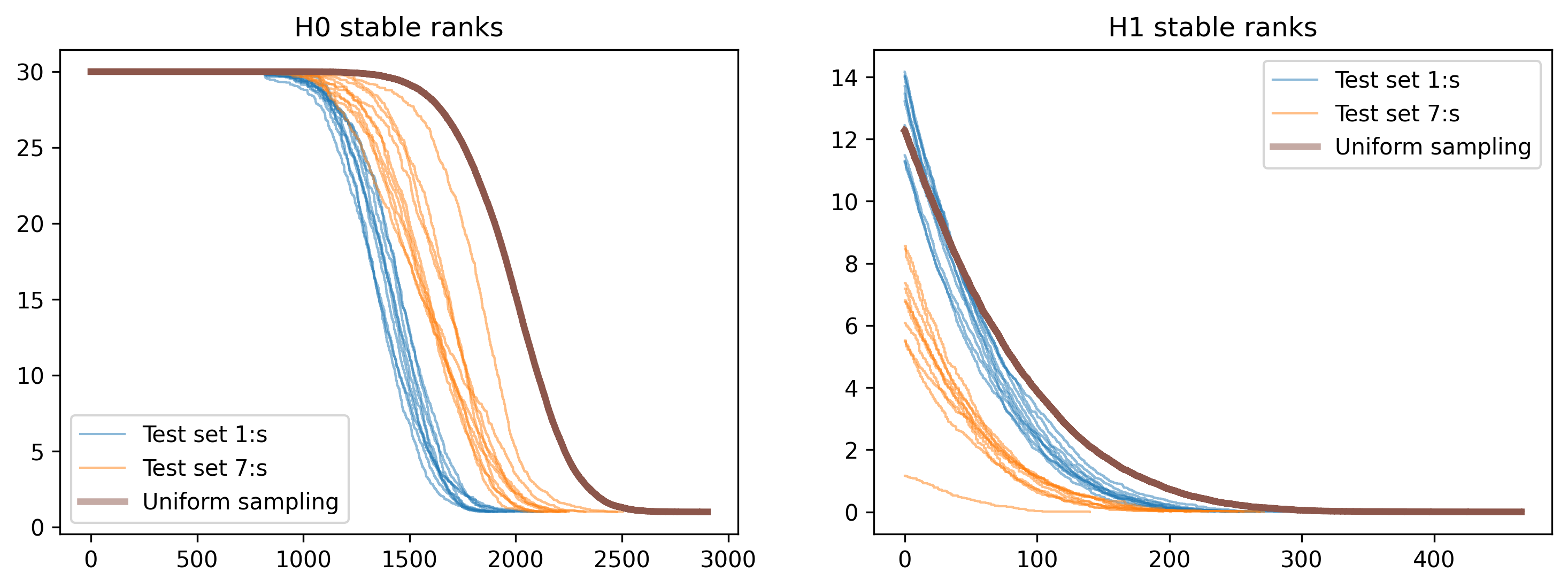}
\end{center}
\caption{Stable ranks corresponding to test set 1:s, test set 7:s and to uniform sampling, for different homological degrees. Training set 2:s and 3:s used as reference object.}\label{illus:5}
\end{figure}

In the previous section, we considered a reference object which consisted of samples from the same data distributions (handwritten 1:s and 7:s) as the points that we sampled from and tried to discriminate. Now, while still trying to distinguish between test set samples of digits 1 and 7, we instead take as our reference object the union $\text{Train}_2\cup\text{Train}_3$ of all samples from the MNIST training set corresponding to digits 2 or 3. Stable ranks are computed following the same procedure, however, for the {\em homological degree} $1$,  we used $1900$ for the {\em truncation} parameter. The obtained stable ranks are illustrated in Figure \ref{illus:5}. Interestingly, when subsampled from different points representing 1:s and 7:s, the geometry of this reference object, which a priori is not related to the data distribution of those digits, nonetheless contains information about those points.

\section{Discussion}
Extracting stable ranks is a simplifying procedure. Finding appropriate parameters 
controlling stable ranks so that relevant aspects of the problem at hand  are retained is the key challenge. In this paper we indicate that choosing an appropriate reference object and ways of sampling it can be used for this purpose. 
For example in experiment \ref{expmnist:unionoftrainingsets} the reference object is the union of the training sets corresponding to  digits $1$ and $7$ which was shown to be effective for distinguishing between  them.  
While analogous  experiments can be repeated with similar results for several other pairs of digits, some pairs of digits were nonetheless harder to distinguish. In experiment \ref{expmnist:unionoftrainingsets} we could in general see that it was harder to distinguish test set digits from the two classes when the global geometries of the digits (see Section \ref{sec:global}) were similar. But while more difficult, it was still often possible, since by sampling from the perspective of different points one can reveal different local geometric patterns that are specific to the digit. A classifier, when fed with such patterns, can thus still learn to distinguish the digits. Moreover, in experiment \ref{expmnist:otherrefobject}, when sampling a reference object that is not the union of the training sets corresponding to the digits we want to distinguish, we are in a different situation where global geometric similarity of the digits does not necessarily matter. Which reference object to choose is however not obvious. Another possibility is to combine different geometric signatures, e.g. stable ranks obtained by taking the training sets corresponding to digits 1 and 7 as separate reference objects, and computed for different homological degrees and parameters. Such signatures could then be combined in e.g. an ensemble learning scheme. We also emphasize that our method by construction only considers relative geometrical aspects to a point. Another interesting direction is thus to combine it with other methods (distance-based machine learning methods, neural networks, etc.) and analyze the combined effect.

\section{Appendix: hierarchical stabilisation}
\label{sec:adfgdfagsdfg} In this appendix we briefly recall the role the parameters {\em density} and {\em truncation}  of our pipeline (see Section~\ref{sec:pipeline}) play  for constructing stable ranks. We refer to~\cite{10.3389/fams.2021.668046, MR4057607,  OliverWojtek} where more information about stable ranks can be found. 

Stable ranks are built using a process called hierarchical stabilization. An input for this process has two ingredients. One is a discrete invariant such as the  rank function $\text{rank}\colon\text{Tame}([0,\infty), \text{vect}_{\mathbb F})\to\mathbb{N}$, which  assigns
to a persistence module its  minimal number of generators.
The other ingredient is a pseudometric  $d$ on the domain of the discrete invariant, which in the case of the rank function is given by  persistence modules  $\text{Tame}([0,\infty), \text{vect}_{\mathbb F})$. 
The outcome of the 
hierarchical stabilization, for the mentioned rank function,  is  a  Lipschitz  function $\widehat{\text{rank}}_d\colon\text{Tame}([0,\infty),  \text{vect}_{\mathbb F})\to\mathcal{M}$, called {\em stable rank}, where $\mathcal{M}$  is the space of Lebesgue measurable functions $[0,\infty )\to [0,\infty)$.
We think about the stable rank function as the model associated
 to the pseudometric $d$. In this framework (supervised) persistence analysis is about identifying these pseudometrics $d$ for which structural properties of the (training)
data are reflected by the geometry of its image in $\mathcal{M}$ through the function $\widehat{\text{rank}}_d$. 

The reason we care about densities and truncations is because any choice of them leads to a pseudometric on persistence modules. Thus we can use densities and truncations as parameters of a rich space of such pseudometrics.
We refer the reader to the mentioned sources for an explanation of how a density and a truncation leads to a pseudometric.  See~\cite{10.3389/fams.2021.668046, MR4057607} for examples where choosing an appropriate density leads to improvement in certain classifications tasks. In this article we have seen that a choice of truncation can also lead to better results.

\bibliographystyle{alpha}
\bibliography{bibliography}

\newcommand{\etalchar}[1]{$^{#1}$}
\begin{thebibliography}{AEK{\etalchar{+}}17}

\bibitem[AEK{\etalchar{+}}17]{MR3625712}
Henry Adams, Tegan Emerson, Michael Kirby, Rachel Neville, Chris Peterson,
  Patrick Shipman, Sofya Chepushtanova, Eric Hanson, Francis Motta, and Lori
  Ziegelmeier.
\newblock Persistence images: a stable vector representation of persistent
  homology.
\newblock {\em J. Mach. Learn. Res.}, 18:Paper No. 8, 35, 2017.

\bibitem[ARSC21]{10.3389/fams.2021.668046}
Jens Agerberg, Ryan Ramanujam, Martina Scolamiero, and Wojciech Chachólski.
\newblock Supervised learning using homology stable rank kernels.
\newblock {\em Frontiers in Applied Mathematics and Statistics}, 7:39, 2021.

\bibitem[Bub15]{MR3317230}
Peter Bubenik.
\newblock Statistical topological data analysis using persistence landscapes.
\newblock {\em J. Mach. Learn. Res.}, 16:77--102, 2015.

\bibitem[CFL{\etalchar{+}}15]{chazal:hal-01073073}
Fr{\'e}d{\'e}ric Chazal, Brittany~Terese Fasy, Fabrizio Lecci, Bertrand Michel,
  Alessandro Rinaldo, and Larry Wasserman.
\newblock {Subsampling Methods for Persistent Homology}.
\newblock In {\em {International Conference on Machine Learning (ICML 2015)}},
  Lille, France, July 2015.

\bibitem[CR20]{MR4057607}
Wojciech Chach\'{o}lski and Henri Riihim\"{a}ki.
\newblock Metrics and stabilization in one parameter persistence.
\newblock {\em SIAM J. Appl. Algebra Geom.}, 4(1):69--98, 2020.

\bibitem[EH08]{MR2405684}
Herbert Edelsbrunner and John Harer.
\newblock Persistent homology---a survey.
\newblock In {\em Surveys on discrete and computational geometry}, volume 453
  of {\em Contemp. Math.}, pages 257--282. Amer. Math. Soc., Providence, RI,
  2008.

\bibitem[Ghr08]{MR2358377}
Robert Ghrist.
\newblock Barcodes: the persistent topology of data.
\newblock {\em Bull. Amer. Math. Soc. (N.S.)}, 45(1):61--75, 2008.

\bibitem[Hau95]{MR1368659}
Jean-Claude Hausmann.
\newblock On the {V}ietoris-{R}ips complexes and a cohomology theory for metric
  spaces.
\newblock In {\em Prospects in topology ({P}rinceton, {NJ}, 1994)}, volume 138
  of {\em Ann. of Math. Stud.}, pages 175--188. Princeton Univ. Press,
  Princeton, NJ, 1995.

\bibitem[LCB]{MNIST}
Yann LeCun, Corinna Cortes, and Christopher~J.C. Burges.
\newblock {\em The MNIST database of handwritten digits}.

\bibitem[OW17]{OliverWojtek}
Gäfvert Oliver and Chach\'olski Wojciech.
\newblock Stable invariants for multidimensional persistence.
\newblock {\em arXiv:1703.03632}, 2017.

\bibitem[SCL{\etalchar{+}}17]{MR3735858}
Martina Scolamiero, Wojciech Chach\'{o}lski, Anders Lundman, Ryan Ramanujam,
  and Sebastian \"{O}berg.
\newblock Multidimensional persistence and noise.
\newblock {\em Found. Comput. Math.}, 17(6):1367--1406, 2017.

\bibitem[Wei11]{MR2777589}
Shmuel Weinberger.
\newblock What is{$\ldots$}persistent homology?
\newblock {\em Notices Amer. Math. Soc.}, 58(1):36--39, 2011.

\end{thebibliography}

\end{document}